\documentclass[12pt]{amsart}
\usepackage[usenames]{color} 
\usepackage{times}
\usepackage{amsfonts}
\usepackage{amsthm}
\usepackage{amsmath}
\usepackage{psfrag}

\usepackage{times}
\usepackage{latexsym,color}
\usepackage{amssymb}
\usepackage{graphicx}

\usepackage{epsfig}

\usepackage{gastex}

\usepackage{tikz}
\usepackage[all]{xy}
 \usepackage[hang]{caption}
 \usetikzlibrary{snakes}
\usetikzlibrary{arrows}

\advance \topmargin by -\headheight
\advance \topmargin by -\headsep
\evensidemargin \oddsidemargin
\def \ind{_{n \in {\mbox{\rm {\scriptsize I$\!$N}}}}}

\newcommand{\GR}{{\mathbb R}}
\newcommand{\GZ}{{\mathbb Z}}

\newcommand{\GN}{{\mathbb N}}

\newcommand{\gN}{\mbox{\rm \scriptsize I$\!$N}}

\newcommand{\ab}{|}

\newtheorem{theorem}{Theorem}
\newtheorem{lemma}[theorem]{Lemma}
\newtheorem{corollary}[theorem]{Corollary}
\newtheorem{proposition}[theorem]{Proposition}

\newtheorem{definition}{Definition}
\newtheorem{remark}{Remark}

\begin{document}

\title[AR iet]{Arnoux-Rauzy interval exchange transformations}

\author[P. Arnoux]{Pierre Arnoux} 
\address{Aix Marseille Universit\'e, CNRS, Centrale Marseille, Institut de Math\' ematiques de Marseille, I2M - UMR 7373\\13453 Marseille, France.}
\email{pierre@pierrearnoux.fr}
\author[J. Cassaigne]{Julien Cassaigne} 
\address{Aix Marseille Universit\'e, CNRS, Centrale Marseille, Institut de Math\' ematiques de Marseille, I2M - UMR 7373\\13453 Marseille, France.}
\email{julien.cassaigne@math.cnrs.fr}
\author[S. Ferenczi]{S\'ebastien Ferenczi} 
\address{Aix Marseille Universit\'e, CNRS, Centrale Marseille, Institut de Math\' ematiques de Marseille, I2M - UMR 7373\\13453 Marseille, France.}
\email{ssferenczi@gmail.com}
\author[P. Hubert]{Pascal Hubert} 
\address{Aix Marseille Universit\'e, CNRS, Centrale Marseille, Institut de Math\' ematiques de Marseille, I2M - UMR 7373\\13453 Marseille, France.}
\email{hubert.pascal@gmail.com}

\date{June 19, 2019}

\begin{abstract} The Arnoux-Rauzy systems are defined in \cite{ar}, both as symbolic systems on three letters and exchanges of six intervals on the circle. In connection with a conjecture of S.P. Novikov, we  investigate the dynamical properties of the interval exchanges, and precise their relation with the symbolic systems, which was known only to be a semi-conjugacy; in order to do this, we define a new system which is an exchange of nine intervals on the line (it was described in \cite{abb} for a particular case). Our main result is that the semi-conjugacy determines a measure-theoretic isomorphism (between the three systems) under a diophantine (sufficient) condition, which is satisfied by almost all Arnoux-Rauzy systems for a  suitable measure; but, under another condition, the interval exchanges are not uniquely ergodic and the isomorphism does not hold for all invariant measures; finally, we give conditions for these  interval exchanges to be weakly mixing.
\end{abstract}
\maketitle

{\em Arnoux-Rauzy dynamical systems} were introduced in \cite{ar} in order to generalize the fruitful triple interaction between Sturmian sequences and rotation of the $1$-torus through the Euclid continued fraction approximation. Arnoux-Rauzy sequences are defined through word-combinatorial conditions, see Section \ref{sar3} below, and have been studied from the combinatorial point of view by many authors, see for example \cite{bulu} \cite{daz} \cite{del} \cite{jupi},  and many others. These sequences  constitute also the restriction to three-letter alphabets of the class of  {\em Episturmian} sequences, defined in \cite{djp}, and extensively  studied, see the surveys \cite{be} and \cite{gjp}.

The first and foremost question on Arnoux-Rauzy sequences was to get  a geometric representation of the associated symbolic dynamical system, the preferred one being as a natural coding of a rotation of the $2$-torus.  The set of possible angles of this rotation is known as the {\em Rauzy gasket}, and defined in Section \ref{sar6} below. A famous particular case, the {\em Tribonacci} sequence, was shown in \cite{rau} to be a {\em natural coding  of a rotation of the $2$-torus}, and thus the corresponding  system is   measure-theoretically isomorphic to that rotation. This was generalized to a larger class of Arnoux-Rauzy systems in \cite{abi}, and recently to almost all Arnoux-Rauzy systems \cite{bst}, in the sense of Definition \ref{aa} below. On the other hand, \cite{cfme} provides counter-examples where this isomorphism cannot hold, see Section \ref{wm} below. For a general Arnoux-Rauzy system, one has to be content with what looks like a second-best geometric representation built in \cite{ar}, a coding of a six-interval exchange on the circle, see Section \ref{sar6} below; in the Tribonacci case, this is the {\em Arnoux-Yoccoz interval exchange} \cite{arn}, linked with the pseudo-Anosov map defined in \cite{ay}. Note that it is still an open question to find other  geometric models, in particular for those Arnoux-Rauzy systems which are not natural codings of rotations of the $2$-torus,  see for example \cite{sir}. 

However, these six-interval exchanges have been recently understood to represent by themselves a very interesting family of systems, as  the dimension over the rationals of the set of lengths of the intervals is quite smaller than the number of intervals (namely, three versus six). This kind of interval exchanges was pointed out (in a very different context and language) by S.P. Novikov, see  \cite{nov}\cite{dn}\cite{mano}, also \cite{d1}\cite{d2}\cite{dedi}. This  prompted several authors to make deep studies of the Rauzy gasket in \cite{le} (Lemma 5.9, attributed to J.-C. Yoccoz) \cite{as} \cite{ahs0} \cite{ahs1}  \cite{gm},  partially solving a conjecture of Novikov, and to look at everything that can be found about this particular family. But indeed, a priori not much is known, as these six-interval exchanges (called AR6 in the present paper) are only semi-conjugate to the original Arnoux-Rauzy systems (called AR3 in the present paper): namely, an AR6 interval exchange admits a coding by a partition into three sets which is an AR3 symbolic system, but this partition is not necessarily a generating partition, while, as far as we know,  the coding by the natural partition  into six intervals of the circle cannot be built by {\em substitutions}, contrarily to its AR3 coding. Hence no property of an AR6 interval exchange can be directly carried out from the underlying AR3 symbolic system. Moreover, while all AR3 systems are known to be  {\em minimal} \cite{ar} and {\em uniquely ergodic} (by Boshernitzan's result  \cite{bos2} using the fact that the language complexity is $2n+1$), in stark contrast, deep geometric methods have allowed I.A. Dynnikov and A. Skripchenko \cite{ds} to prove, again in a completely different language, the existence of minimal  non-uniquely ergodic AR6 interval exchanges.\\

The relation between AR6 interval exchanges and underlying AR3 symbolic systems was partially tackled in \cite{abb}, though only in the particular  case of Tribonacci, and with a certain lack of details: that paper defines yet another Arnoux-Rauzy interval exchange, this time on nine intervals (called AR9 in the present paper), where an AR3 appears again as a coding by  a partition into three sets, and where the coding by the natural partition into nine intervals can be explicitly generated by  a substitution. This is the key for studying ergodic properties of AR9 interval exchanges, and extending them to the AR6 interval exchanges which appear as factors of AR9. The one stated in \cite{abb} is the measure-theoretic isomorphism between the three corresponding systems (AR3, AR6, AR9) in the Tribonacci case, though no proof is offered. \\ 

In the present paper, we generalize the construction of AR9 systems to every set of parameters in the Rauzy gasket, and their construction by substitutions, using an {\em induction} process defined in Section \ref{sind} below; we use them to derive dynamical properties of AR6 and AR9 systems. Our main result is 

\begin{theorem}\label{main} Almost every (in the sense of \cite{bst}, see Definition \ref{aa} below) AR9 or AR6 interval exchange is uniquely ergodic and measure-theoretically isomorphic to its AR3 coding. \end{theorem}

This theorem could be deduced (using Lemma \ref{tourc} below and some extra work) from the ergodicity of the induction process;  we choose to derive it from a stronger result, namely an explicit sufficient  diophantine condition (Proposition \ref{mtours} and Theorem \ref{imp} below) for measure-theoretic isomorphism between the corresponding AR9, AR6 and AR3 systems, which also implies unique ergodicity for the AR6 and AR9. This condition is satisfied by almost all Arnoux-Rauzy systems (Proposition \ref{pt} below), and many explicit examples including Tribonacci, all systems which are periodic points under the induction, and, more generally,  all the so-called Arnoux-Rauzy systems with {\em bounded weak partial quotients} (Proposition \ref{bqp} below). Thus

\begin{corollary}\label{cp1} Almost all AR9 or AR6 interval exchanges, including the Tribonacci ones and all those with bounded  weak partial quotients, are measure-theoretically isomorphic to rotations of the $2$-torus. \end{corollary}

Thus at last we have proved the isomorphism result for the  Tribonacci case; this  provides the backbone of an answer to Question 9 (asked by G. Forni) in \cite{fh}  and this was another motivation for the present paper.

\begin{corollary}\label{cp2} The Arnoux-Yoccoz interval exchange, or else the Tribonacci AR9,  provide nontrivial examples of {\em rigid}  (a sequence of powers converges to the identity in $\mathcal L^2$) self-induced interval exchanges. \end{corollary}

 Then we give a class of examples of non-uniquely ergodic AR9 (or AR6) which may be somewhat more explicit than those in \cite{ds}, and give both examples and counter-examples to the isomorphism problem: these AR9 are measure-theoretically isomorphic to their AR3 coding if we equip them with an ergodic invariant measure, but of course this cannot hold if we take one of the many non-ergodic measures. Then we show that weak mixing is also present in the class of AR9 (or AR6) systems.\\

{\bf Acknowledgement:} this research was born from a discussion with V. Delecroix during the FWF/JSPS project meeting in Salzburg in 2018; a part of it was carried out when the second and third authors participated in the meeting organized by S. Brlek in Murter (Croatia) in april 2018, and another part in july 2018 while the first author was in Unit\' e Mixte IMPA-CNRS (Institut Jean-Christophe Yoccoz) in Rio de Janeiro and
the third author was a temporary visitor of IMPA through the R\' eseau Franco-Br\' esilien en Math\' ematiques.\\

\section{Basic definitions}
We look at finite {\em words} on a finite alphabet ${\mathcal A}=\{1,...k\}$. A word $w_1...w_t$ has
{\em length} $\ab w\ab=t$. The {\em concatenation} of two words $w$ and $w'$ is denoted by $ww'$. 

\begin{definition}\label{dln} 
\noindent A word $w=w_1...w_t$ {\em occurs at place $i$} in a word $v=v_1...v_s$ or an infinite sequence  $v=v_1v_2...$ if $w_1=v_i$, ...$w_t=v_{i+t-1}$. We say that $w$ is a {\em subword} of $v$. \end{definition}

\begin{definition}
A {\em   language} $L$ over $\mathcal A$ is a set of words. In the present paper, all languages are assumed to be {\em factorial} (if $w$ is in $L$, all its subwords are in $L$), and {\em extendable} (if $w$ is in $L$, $aw$ is in $L$ for at least one letter $a$ of $\mathcal A$, and $wb$ is in $L$ for at least one letter $b$ of $\mathcal A$).

\noindent A language $L$ is
{\em  minimal} if for each $w$ in $L$ there exists $n$ such that $w$
occurs in each word  of $L$ of length $n$.

\noindent The language $L(u)$ of an infinite sequence $u$ is the set of its finite subwords.
\end{definition}

\begin{definition}\label{subs} A {\em substitution} $\psi$ is an application from
an alphabet $\mathcal A$ into the set ${\mathcal A}^{\star}$ of finite words on 
$\mathcal A$; it
extends naturally to a morphism of ${\mathcal A}^{\star}$ for the operation of concatenation. 
\end{definition}

\begin{definition} The {\em symbolic dynamical system} associated to a language $L$ is the one-sided shift $S(x_0x_1x_2...)=x_1x_2...$ on the subset $Y_L$ of ${\mathcal A}^{\GN}$ made with the infinite sequences such that for every $t<s$, $x_t...x_s$ is in $L$.
\end{definition}

Note that the symbolic dynamical system $(X_L,S)$ is {\em minimal} (in the usual sense, every orbit is dense) if and only if the language $L$ is mimimal.

\begin{definition}\label{sy} For a dynamical system $(X',U)$ and a finite partition $\{P_1,\ldots P_l\}$ of $X'$, the {\em trajectory} of a point $x$ in $X'$ is the infinite
sequence
$(x_{n})\ind$ defined by $x_{n}=i$ if $U^nx$ falls into
$P_i$, $1\leq i\leq l$.

\noindent Then if $L$ is the language made of all the finite subwords of all the  trajectories, $(Y_L,S)$ is called the {\em coding} of $(X',U)$ by the partition $\{P_1,\ldots P_l\}$. 
\end{definition}

\section{Classical Arnoux-Rauzy systems}
\subsection{AR3 symbolic systems}\label{sar3}
These systems are the ``genuine" Arnoux-Rauzy systems; we take here as a definition  their constructive characterization, derived in \cite{ar} from the original definition, and modified in the present paper by a renaming of letters and words.

\begin{definition}\label{dar3} An {\em AR3 symbolic system} is the  symbolic system on $\{a,b,c\}$ generated by the three substitutions 
\begin{itemize}
\item $\sigma_I$: $a\to ab$, $b\to ac$, $c\to a$,
\item $\sigma_{II}$: $a\to ab$, $b\to a$, $c\to ac$,
\item $\sigma_{III}$: $a\to a$, $b\to ab$, $c\to ac$,
\end{itemize}
and a directing sequence $r_n$,  $n\in \GN^{\star}$, $r_n \in\{I,II,III\}$, taking the value $I$ infinitely many times. 

 Namely, it is  the symbolic system $(Y_3,S)$ whose language is generated by the words $A_k=\sigma_{r_1}...\sigma_{r_k}a$, $B_k=\sigma_{r_1}...\sigma_{r_k}b$, $C_k=\sigma_{r_1}...\sigma_{r_k}c$, $k\geq 1$. The respective lengths of the words $A_k$, $B_k$, $C_k$ will always be denoted by $h_{a,k}$, $h_{b,k}$, $h_{c,k}$.
\end{definition}

As mentioned in the introduction,  $(Y_3,S)$ is minimal, and {\em uniquely ergodic}: there is a  unique invariant probability measure, denoted by  $\mu$.\\

Note that our modification  of the rules changes the usual condition of \cite{ar}, that each substitution is used infinitely often, to the present condition that $\sigma_I$ is used infinitely often.
The most famous particular case is the {\em Tribonacci system}, where $r_n=I$ for all $n$.

 \subsection{Partial quotients and multiplicative rules}\label{pqm}
These quantities are defined in \cite{cfme}, but we redefine them here as the notations are different.

\begin{definition}\label{d2} We write the directing sequence $(r_n)$  in a unique way  as $k_1-1\geq 0$ times the symbol $III$ followed by one symbol $I$ or $II$, then $k_2-1\geq 0$ times  $III$ followed by one  $I$ or $II$ etc.... the $k_n\geq 1$ are then called the {\em partial quotients} of the system.\\
The {\em multiplicative times} are $m_0=0$, $m_n=k_1+...k_n$, $n\geq 1$: they are the times $m$ for which $r_m\neq III$.
\end{definition}

Then the words $A_{m_n}$, $B_{m_n}$, $C_{m_n}$ can be built by the following {\em multiplicative rules}, which could also be expressed by substitutions but would need a countable set of them:

\begin{itemize} \item if  $r_{m_{n+1}}=I$, we say
 that {\em the $n+1$-th multiplicative rule is a rule  $I_m$}, and we have
\begin{itemize} \item  $A_{m_{n+1}}=A_{m_n}^{k_{n+1}}B_{m_n}$, \item $B_{m_{n+1}}=A_{m_n}^{k_{n+1}}C_{m_n}$, \item $C_{m_{n+1}}=A_{m_n}$; \end{itemize} 
 \item if  $r_{m_{n+1}}=II$, we say
 that {\em the $n+1$-th multiplicative rule is a rule  $II_m$}, and
\begin{itemize} \item
 $A_{m_{n+1}}=A_{m_n}^{k_{n+1}}B_{m_n}$, \item $B_{m_{n+1}}=A_{m_n}$, \item $C_{m_{n+1}}=A_{m_n}^{k_{n+1}}C_{m_n}$.\end{itemize}
 \end{itemize}
 
   For Tribonacci, we have $k_n=1$ for all $n$, and all multiplicative rules are  $I_m$.

We recall that in \cite{cfme}, we  use different substitutions (called ``(additive) concatenation rules" in that paper), and the sequence of multiplicative rules (as defined in that paper) corresponds to the successive number of times we use each substitution: the $n+1$-th multiplicative rule is of type $1$ whenever
the $m_{n-1}$-th and $m_{n+1}$-th substitutions are different. Then the
 $H_n$, $G_n$ and $J_n$ of \cite{cfme} are exactly the same as  respectively $A_{m_n}$,  $B_{m_n}$ and $C_{m_n}$ in the present paper, and types 1 and 2 of \cite{cfme} correspond to our rules $I_m$ and $II_m$.

We shall use the inequalities proved in Lemma 7 of \cite{cfme} at the multiplicative times: namely $h_{b,m_n}\leq 2h_{a,m_n}$ and $h_{c,m_n}\leq 2h_{a,m_n}$. These are not true in general at other (additive) times $p\neq m_n$.

 \subsection{AR6 interval exchanges}\label{sar6}

These exchanges of six intervals on a circle are defined in \cite{arn} for Tribonacci, see also \cite{ay}, and \cite{ar} for the general case. 

\begin{definition} The {\em Rauzy gasket} $\Gamma$ is the set of triples of positive real 
  numbers $(a_0,b_0,c_0)$, such that, if we define recursively the numbers $a_n$, $b_n$, $c_n$ by taking the triple $(a_{n-1}-b_{n-1}-c_{n-1},b_{n-1}, c_{n-1})$ and reordering it, then for each $n\geq 0$ we have $a_n>b_n>c_n>0$. \end{definition}
 
\begin{definition}\label{dar6} An {\em AR6 interval exchange}  $(X_6,T)$
is defined in the following way from any triple $(a_0,b_0,c_0)$ in $\Gamma$:
$X_6$ is a circle of length $2a_0+2b_0+2c_0$. The circle is partitioned into three intervals of respective lengths $2a_0$, $2b_0$, $2c_0$, then each one is cut into two halves; the action of $T$ first exchanges by translations respectively the two intervals of length $a_0$, the two intervals of length $b_0$, the two intervals of length $c_0$, then translates everything by $a_0+b_0+c._0$, i.e. a half-circle.\end{definition}

We could also look at the same transformation as an exchange of seven intervals on the interval $[0,2a_0+2b_0+2c_0[$, but a better model on the interval will be given in Section \ref{sar9} below. 

Note that he location of the origin on the circle does not change the system up to topological conjugacy and  measure-theoretic isomorphism for any invariant measure (in the sense that any invariant measure on one of them can be carried to the other one, and the two measure-theoretic systems are isomorphic). Similarly, the order between  the intervals of lengths $2a_0$, $2b_0$, $2c_0$ on the circle is not mentioned in Definition \ref{dar6} (the fact that it is not always the same is somewhat understated  in \cite{ar}); by changing the origin, we can reduce the number of possible orders to two, and the two  AR6 interval exchanges defined with the same  $(a_0,b_0,c_0)$ but different orders of these intervals are conjugate by a symmetry on the circle, thus are also topologically conjugate and measure-theoretically isomorphic for any invariant measure.\\

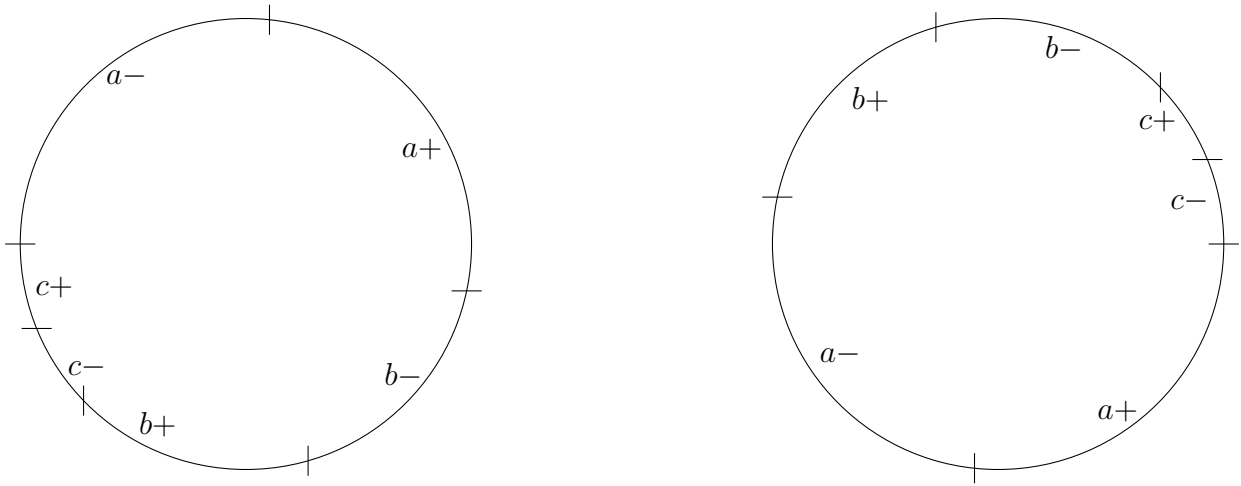
\begin{figure}[h] \label{fig: AR6 interval exchange}
\begin{center}
\begin{tikzpicture}[scale = 1]

\draw (-5,0) circle(3);

\draw({-5.2-3*cos(0)}, {3*sin(0)})-- ({-4.8-3*cos(0)}, {3*sin(0)}); 
\draw({-5-3*cos(96)}, {-.2+3*sin(96)})-- ({-5-3*cos(96)}, {.2+3*sin(96)}); 
\draw({-5.2-3*cos(192)}, {3*sin(192)})-- ({-4.8-3*cos(192)}, {3*sin(192)}); 
\draw({-5-3*cos(254)}, {-.2+3*sin(254)})-- ({-5-3*cos(254)}, {.2+3*sin(254)}); 
\draw({-5-3*cos(316)}, {-.2+3*sin(316)})-- ({-5-3*cos(316)}, {.2+3*sin(316)}); 
\draw({-5.2-3*cos(338)}, {3*sin(338)})-- ({-4.8-3*cos(338)}, {3*sin(338)});

\draw({-5-3*cos(48)}, {3*sin(48)}) node[above, right]{$a-$};
\draw({-5.5-3*cos(144)}, {-.5+3*sin(144)}) node[above, right]{$a+$};
\draw({-4.7-3*cos(223)}, {.3+3*sin(223)}) node[above, left]{$b-$};
\draw({-5-3*cos(285)}, {.5+3*sin(285)}) node[above, left]{$b+$};
\draw({-5-3*cos(327)}, {3*sin(327)}) node[above, right]{$c-$};
\draw({-5-3*cos(349)}, {3*sin(349)}) node[above, right]{$c+$};

 \draw (5,0) circle(3);
 \draw({5.2-3*cos(180)}, {3*sin(180)})-- ({4.8-3*cos(180)}, {3*sin(180)}); 
\draw({5-3*cos(276)}, {-.2+3*sin(276)})-- ({5-3*cos(276)}, {.2+3*sin(276)}); 
\draw({5.2-3*cos(12)}, {3*sin(12)})-- ({4.8-3*cos(12)}, {3*sin(12)}); 
\draw({5-3*cos(74)}, {-.2+3*sin(74)})-- ({5-3*cos(74)}, {.2+3*sin(74)}); 
\draw({5-3*cos(136)}, {-.2+3*sin(136)})-- ({5-3*cos(136)}, {.2+3*sin(136)}); 
\draw({5.2-3*cos(158)}, {3*sin(158)})-- ({4.8-3*cos(158)}, {3*sin(158)});

\draw({5-3*cos(228)}, {3*sin(228)}) node[above, left]{$a+$};
\draw({5.2-3*cos(335)}, {-.2+3*sin(335)}) node[above, right]{$a-$};
\draw({6-3*cos(40)}, {3*sin(40)}) node[below, left]{$b+$};
\draw({5.5-3*cos(105)}, {-.3+3*sin(105)}) node[below, left]{$b-$};
\draw({5-3*cos(147)}, {3*sin(147)}) node[above, left]{$c+$};
\draw({5-3*cos(169)}, {3*sin(169)}) node[above, left]{$c-$};

\end{tikzpicture}

\caption{AR6 interval exchange}
\end{center}

\end{figure}
 
 For example, when the intervals  of lengths $2a_0$, $2b_0$, $2c_0$ are successive intervals of the circle in that order, $T$ is shown in Figure 1, where on the left circle  $a-$, $a+$, $b-$, ... denote the intervals of length $a_0$, $a_0$, $b_0$ ... while on the right circle the letters correspond to the images of these intervals by the transformation.
 If in Figure 1 we choose to put the origin at the left end of the interval denoted by $a-$,
 $[0,a_0)$ is sent to $[a_0+a_0+b_0+c_0, 2a_0+a_0+b_0+c_0)$ modulo $2a_0+2b_0+2c_0$,  
 $[a_0,2a_0)$ is sent to $[a_0+b_0+c_0, a_0+a_0+b_0+c_0)$ modulo $2a_0+2b_0+2c_0$, etc...  \\
 
 The link between AR3 symbolic systems and AR6 interval exchanges, studied in \cite{ar}, will be described in Section \ref{nno} below. But, as pointed out in the introduction, we do not know any constructive way to build directly the language of
the {\em natural} coding of the AR6 interval exchange, that is its coding by the partition into its six intervals of continuity on the circle, coded by $a-$, $a+$, $b-$, $b+$, $c-$, $c+$. That is why we need to introduce one more class of Arnoux-Rauzy systems.

\subsection{Note on endpoints}\label{nep}
One recurring problem when dealing with interval exchanges is what to do with interval endpoints? A satisfying answer to this question is given by M. Keane in  Section 5 of \cite{kea1}: by carefully doubling the endpoints and their orbits, he defines a Cantor set on which the transformation becomes an homeomorphism, and show this is equivalent to taking the natural coding by the partition into defining intervals.  In the present paper, to make definitions easier,  we do not use Keane's construction, and all intervals are closed on the left, open on the right; but that will introduce technical difficulties, see Remark \ref{rep} below. 

\section{The new systems: Arnoux-Rauzy on nine symbols}\label{sar9}
\subsection{AR9 interval exchanges}

These are defined for the particular case of Tribonacci in \cite{abb}. Here we define them  in full generality, in a deliberately pedestrian way, which does not reveal  how they were devised; the grand geometry underlying and motivating the construction, generalizing the geometry in \cite{abb}, will appear in a further paper. Note that we use the same symbol $T$ for AR9 and AR6 interval exchanges in view of Proposition \ref{ar63} below.\\

An AR9 interval exchange is defined by a point $(a_0,b_0,c_0)$ in $\Gamma$, as an exchange  of nine intervals on a union of three disjoint intervals on the line.

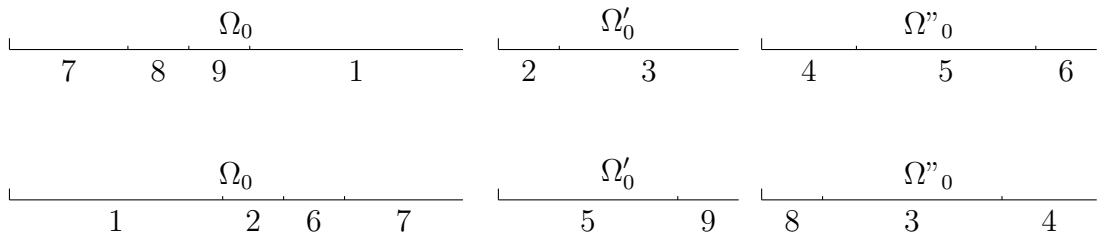
\begin{figure}[h] \label{fig: AR9 interval exchange}
\begin{center}
\begin{tikzpicture}[scale = 5]

\draw(0,.1)--(0,.13);

\draw(0,.5)--(0,.53);
\draw (1.3,.1) node[above]{};
\draw(1.3,.1)--(1.3,.13);
\draw(1.3,.5)--(1.3,.53);
\draw (.603,.5) node[above]{$\Omega_0$};
\draw (.603,.1) node[above]{$\Omega_0$};
\draw (1.619,.5) node[above]{$\Omega'_0$};
\draw (1.619,.1) node[above]{$\Omega'_0$};
\draw (2.445,.5) node[above]{$\Omega"_0$};
\draw (2.445,.1) node[above]{$\Omega"_0$};
\draw (2,.1) node[above]{};
\draw(2,.1)--(2,.13);
 
\draw(2,.5)--(2,.53);

\draw(0,.1)--(1.206,.1);
\draw(.283,.1) node[below]{$1$};
\draw(.567,.1) node[above]{};
\draw(.567,.1)--(.567,.11);
\draw (.648,.1) node[below]{$2$};
\draw(.729,.1) node[above]{};
\draw(.729,.1)--(.729,.11);
\draw (.81,.1) node[below]{$6$};
\draw(.891,.1) node[above]{};
\draw(.891,.1)--(.891,.11);
\draw(1.048,.1) node[below]{$7$};

\draw(0,.5)--(1.206,.5);
\draw (.157,.5) node[below]{$7$};
\draw(.315,.5) node[above]{};
\draw(.315,.5)--(.315,.51);
\draw (.396,.5) node[below]{$8$};
\draw(.477,.5) node[above]{};
\draw(.477,.5)--(.477,.51);
\draw (.558,.5) node[below]{$9$};
\draw(.639,.5) node[above]{};
\draw(.639,.5)--(.639,.51);
\draw(.922,.5) node[below]{$1$};
\draw(.1206,.5) node[above]{};

\draw(1.3,.1)--(1.939,.1);
\draw (1.538,.1) node[below]{$5$};
\draw(1.777,.1) node[above]{};
\draw(1.777,.1)--(1.777,.11);
\draw (1.858,.1) node[below]{$9$};

\draw(1.3,.5)--(1.939,.5);
\draw (1.381,.5) node[below]{$2$};
\draw(1.462,.5) node[above]{};
\draw(1.462,.5)--(1.462,.51);
\draw (1.7,.5) node[below]{$3$};
\draw(1.939,.5) node[above]{};

\draw(2,.5)--(2.891,.5);
\draw (2.126,.5) node[below]{$4$};
\draw(2.252,.5) node[above]{};
\draw(2.252,.5)--(2.252,.51);
\draw (2.49,.5) node[below]{$5$};
\draw(2.729,.5) node[above]{};
\draw(2.729,.5)--(2.729,.51);
\draw (2.81,.5) node[below]{$6$};

\draw(2,.1)--(2.891,.1);
\draw (2.081,.1) node[below]{$8$};
\draw(2.162,.1) node[above]{};
\draw(2.162,.1)--(2.162,.11);
\draw (2.4,.1) node[below]{$3$};
\draw(2.639,.1) node[above]{};
\draw(2.639,.1)--(2.639,.11);
\draw (2.765,.1) node[below]{$4$};

\end{tikzpicture}\\
\caption{AR9 interval exchange}

\end{center}

\end{figure}

\begin{definition}\label{dar9} For a point $(a_0,b_0,c_0)$ in $\Gamma$, an AR9 interval exchange $(X_9,T)$ is defined on the union of 
three disjoint intervals $\Omega_0$ of length $a_0+b_0$, $\Omega'_0$ of length $b_0+c_0$, $\Omega"_0$ of length $a_0+c_0$.

An AR9 interval exchange in the first, second or third order is defined in the following way:
\begin{itemize}
\item in the first order, from left to right we see $\Omega_0$, $\Omega'_0$, $\Omega"_0$;  in the second order, from left to right we see $\Omega'_0$, $\Omega"_0$, $\Omega_0$;  in the third order, from left to right we see $\Omega"_0$, $\Omega_0$, $\Omega'_0$,
\item we partition the interval  $\Omega_0$,  from left to right, into four intervals of successive lengths $b_0-c_0$, $c_0$, $c_0$, $a_0-c_0$, denoted respectively by $I_{7,0}$, $I_{8,0}$, $I_{9,0}$, $I_{1,0}$, and into four intervals of successive lengths $a_0-c_0$, $c_0$, $c_0$, $b_0-c_0$, which we define respectively to be $TI_{1,0}$, $TI_{2,0}$, $TI_{6,0}$, $TI_{7,0}$,
\item we partition the interval $\Omega'_0$, from left to right, into two intervals of successive lengths $c_0$, $b_0$,  denoted respectively by $I_{2,0}$, $I_{3,0}$, and into two intervals of successive lengths $b_0$, $c_0$,  which we define respectively to be $TI_{5,0}$, $TI_{9,0}$,
\item we partition the interval $\Omega"_0$ from left to right, into three intervals of successive lengths $a_0-b_0$, $b_0$, $c_0$, denoted respectively by $I_{4,0}$, $I_{5,0}$, $I_{6,0}$, and into three intervals of successive lengths $c_0$, $b_0$, $a_0-b_0$, which we define respectively to be $TI_{8,0}$, $TI_{3,0}$, $TI_{4,0}$.
\end{itemize} 

An AR9 interval exchange in the reversed first, second or third order is defined in the same way, except that in all items above``from left to right" is replaced by ``from right to left" (note that all intervals are still closed on the left, open on the right).
\end{definition}

It is clear from the definition that two AR9 interval exchanges defined with the same $(a_0,b_0,c_0)$ but different actual locations on the line of the intervals $\Omega_0$, $\Omega'_0$, $\Omega"_0$ (equivalently, different locations of the origin and gaps between the intervals), or different orders, are conjugate by a map which is continuous except on a finite number of points, and measure-theoretically isomorphic for any invariant measure, in the sense of Section \ref{sar6} above; all will be topologically isomorphic  if we suppose  no two of the intervals $\Omega_0$, $\Omega'_0$, $\Omega"_0$ are adjacent. 

We could also define AR9 interval exchanges on the circle, gluing $\Omega_0$, $\Omega'_0$, $\Omega"_0$ as in Proposition \ref{ar63} below, but we prefer to define them on the line as (contrarily to the AR6 case) there is no need to add an interval. If we choose the $\Omega_0$, $\Omega'_0$, $\Omega"_0$ to be adjacent, and this  is allowed by our definition, we get examples of  ``usual" nine-interval exchanges as in \cite{kea1},  defined on one interval; but as we shall see below this adjacency will not be conserved by induction, so we have to use the more general family. We shall check that all our results, in particular Lemma \ref{toursint} below, which states the adjacency of certain intervals, is true whatever the gaps between $\Omega_0$, $\Omega'_0$, $\Omega"_0$.

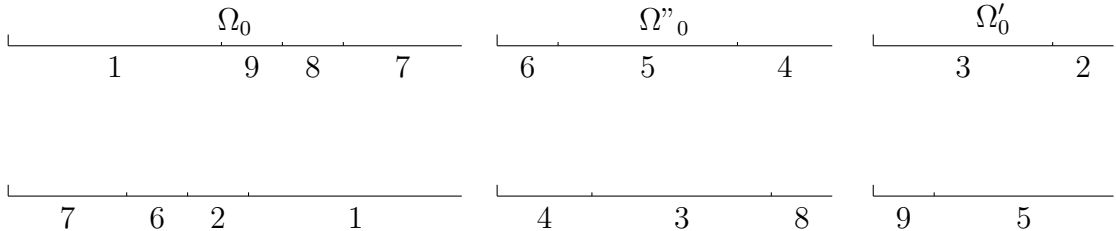
\begin{figure}[h] \label{far9r}
\begin{center}
\begin{tikzpicture}[scale = 5]

\draw (.603,.5) node[above]{$\Omega_0$};

\draw (2.619,.5) node[above]{$\Omega'_0$};

\draw (1.745,.5) node[above]{$\Omega"_0$};

\draw(0,.5)--(0,.53);
\draw(0,.1)--(0,.13);

\draw(1.3,.5)--(1.3,.53);
\draw(1.3,.1)--(1.3,.13);

\draw(2.3,.5)--(2.3,.53);
\draw(2.3,.1)--(2.3,.13);

\draw(0,.1)--(1.206,.1);
\draw(.924,.1) node[below]{$1$};

\draw(.639,.1)--(.639,.11);
\draw (.558,.1) node[below]{$2$};
\draw(.477,.1)--(.477,.11);
\draw (.396,.1) node[below]{$6$};

\draw(.315,.1)--(.315,.11);
\draw(.158,.1) node[below]{$7$};

\draw(2.3,.1)--(2.939,.1);

\draw (2.701,.1) node[below]{$5$};

\draw(2.462,.1)--(2.462,.11);
\draw (2.381,.1) node[below]{$9$};

\draw(1.3,.1)--(2.191,.1);

\draw (2.11,.1) node[below]{$8$};

\draw(2.029,.1)--(2.029,.11);
\draw (1.791,.1) node[below]{$3$};

\draw(1.552,.1)--(1.552,.11);
\draw (1.426,.1) node[below]{$4$};

\draw(0,.5)--(1.206,.5);
\draw (1.049,.5) node[below]{$7$};

\draw(.891,.5) node[above]{};
\draw(.891,.5)--(.891,.51);
\draw (.81,.5) node[below]{$8$};

\draw(.729,.5) node[above]{};
\draw(.729,.5)--(.729,.51);
\draw (.648,.5) node[below]{$9$};

\draw(.567,.5) node[above]{};
\draw(.567,.5)--(.567,.51);
\draw(.284,.5) node[below]{$1$};

\draw(2.3,.5)--(2.939,.5);
\draw (2.858,.5) node[below]{$2$};

\draw(2.777,.5) node[above]{};
\draw(2.777,.5)--(2.777,.51);
\draw (2.539,.5) node[below]{$3$};

\draw(1.939,.5) node[above]{};

\draw(1.3,.5)--(2.191,.5);
\draw (2.065,.5) node[below]{$4$};

\draw(1.939,.5) node[above]{};
\draw(1.939,.5)--(1.939,.51);
\draw (1.701,.5) node[below]{$5$};

\draw(1.462,.5) node[above]{};
\draw(1.462,.5)--(1.462,.51);
\draw (1.381,.5) node[below]{$6$};

\draw(2.891,.5) node[above]{};

\end{tikzpicture}\\
\caption{AR9 interval exchange in reversed order}

\end{center}

\end{figure}

For example, an AR9 interval exchange in the first order is shown  in Figure 2, where $i$ in the upper part corresponds to $I_{i,0}$ and $i$ in the lower part corresponds to $TI_{i,0}$. An example in the reversed second order is shown in Figure 3.

 \subsection{Induction}\label{sind}

Now, we take an AR9 system $(X_9,T)$; to fix ideas, we suppose it is in the first order. Let $T_1$ be the induced map of $T$ on 
$I_{1,0}\cup I_{2,0}\cup I_{3,0}\cup I_{4,0}$. We define $a_1>b_1>c_1$ as the triple $(a_0-b_0-c_0,b_0, c_0)$ after reordering. Then there are three cases, which we tackle by growing  order of difficulty.

\subsubsection{Induction step case III: $a_1=a_0-b_0-c_0.$} Then $b_1=b_0$, $c_1=c_0$.

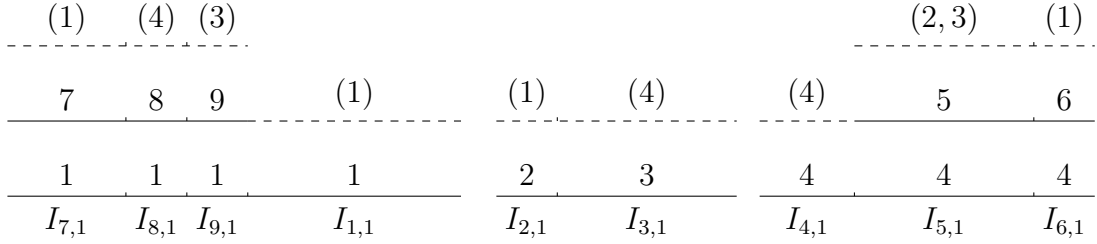
\begin{figure}[h] \label{Induction Case III}
\begin{center}
\begin{tikzpicture}[scale = 5]

\draw (0,.1) node[above]{};
\draw (1.35,.1) node[above]{};
\draw (2.05,.1) node[above]{};

\draw(0,.1)--(1.206,.1);
\draw (.157,.1) node[above]{$1$};
\draw (.157,.1) node[below]{$I_{7,1}$};
\draw(.315,.1) node[above]{};
\draw(.315,.1)--(.315,.11);
\draw (.396,.1) node[above]{$1$};
\draw (.396,.1) node[below]{$I_{8,1}$};
\draw(.477,.1) node[above]{};
\draw(.477,.1)--(.477,.11);
\draw (.558,.1) node[above]{$1$};
\draw (.558,.1) node[below]{$I_{9,1}$};
\draw(.639,.1) node[above]{};
\draw(.639,.1)--(.639,.11);
\draw(.922,.1) node[above]{$1$};
\draw(.922,.1) node[below]{$I_{1,1}$};

\draw(1.3,.1)--(1.939,.1);
\draw (1.381,.1) node[above]{$2$};
\draw (1.381,.1) node[below]{$I_{2,1}$};
\draw(1.462,.1) node[above]{};
\draw(1.462,.1)--(1.462,.11);
\draw (1.7,.1) node[above]{$3$};
\draw (1.7,.1) node[below]{$I_{3,1}$};
\draw(1.939,.1) node[above]{};

\draw(2,.1)--(2.891,.1);
\draw (2.126,.1) node[above]{$4$};
\draw (2.126,.1) node[below]{$I_{4,1}$};
\draw(2.252,.1) node[above]{};
\draw(2.252,.1)--(2.252,.11);
\draw (2.49,.1) node[above]{$4$};
\draw (2.49,.1) node[below]{$I_{5,1}$};
\draw(2.729,.1) node[above]{};
\draw(2.729,.1)--(2.729,.11);
\draw (2.81,.1) node[above]{$4$};
\draw (2.81,.1) node[below]{$I_{6,1}$};
\draw(2.891,.1) node[above]{};

\draw (0,.3) node[above]{};

\draw (1.35,.3) node[above]{};
\draw (2.05,.3) node[above]{};

\draw(0,.3)--(.639,.3);
\draw[dashed](.639,.3)--(1.206,.3);
\draw (.157,.3) node[above]{$7$};
\draw(.315,.3) node[above]{};
\draw(.315,.3)--(.315,.31);
\draw (.396,.3) node[above]{$8$};
\draw(.477,.3) node[above]{};
\draw(.477,.3)--(.477,.31);
\draw (.558,.3) node[above]{$9$};
\draw(.639,.3) node[above]{};
\draw(.922,.3) node[above]{$(1)$};

\draw[dashed](1.3,.3)--(1.939,.3);
\draw (1.381,.3) node[above]{$(1)$};
\draw(1.462,.3) node[above]{};
\draw(1.462,.3)--(1.462,.31);
\draw (1.7,.3) node[above]{$(4)$};
\draw(1.939,.3) node[above]{};

\draw[dashed](2,.3)--(2.252,.3);
\draw(2.252,.3)--(2.891,.3);
\draw (2.126,.3) node[above]{$(4)$};
\draw(2.252,.3) node[above]{};
\draw (2.49,.3) node[above]{$5$};
\draw(2.729,.3) node[above]{};
\draw(2.729,.3)--(2.729,.31);
\draw (2.81,.3) node[above]{$6$};
\draw(2.891,.3) node[above]{};

\draw (0,.5) node[above]{};

\draw (1.35,.5) node[above]{};
\draw (2.05,.5) node[above]{};

\draw[dashed](0,.5)--(.639,.5);
\draw (.157,.5) node[above]{$(1)$};
\draw(.315,.5) node[above]{};
\draw(.315,.5)--(.315,.51);
\draw (.396,.5) node[above]{$(4)$};
\draw(.477,.5) node[above]{};
\draw(.477,.5)--(.477,.51);
\draw (.558,.5) node[above]{$(3)$};
\draw(.639,.5) node[above]{};
\draw(.922,.5) node[below]{};

\draw (1.381,.5) node[below]{};
\draw(1.462,.5) node[above]{};
\draw (1.7,.5) node[below]{};
\draw(1.939,.5) node[above]{};

\draw[dashed](2.252,.5)--(2.891,.5);
\draw (2.126,.5) node[below]{};
\draw(2.252,.5) node[above]{};
\draw (2.49,.5) node[above]{$(2,3)$};
\draw(2.729,.5) node[above]{};
\draw(2.729,.5)--(2.729,.51);
\draw (2.81,.5) node[above]{$(1)$};
\draw(2.891,.5) node[above]{};

\end{tikzpicture}\\
\caption{Induction Case III}

\end{center}

\end{figure}
The situation is essentially described in Figure 4. The induction set $I_{1,0}\cup I_{2,0}\cup I_{3,0}\cup I_{4,0}$ is a disjoint union of three intervals which we denote by $\Omega_1$, $\Omega'_1$, $\Omega"_1$, and is further cut into nine new intervals $I_{i,1}$, whose respective lengths are, from left to right,  $b_1-c_1$, $c_1$, $c_1$, $a_1-c_1$, $c_1$, $b_1$, $a_1-b_1$, $b_1$, $c_1$. Then $T$ acts on the picture as a move upwards, until we reach again the induction set, which is marked by dashed lines.
Each interval of the picture is labelled by  $j$ above if it is in $I_{j,0}$; the labels  are between parentheses for the dashed intervals, as they will not be used further (note that $T_1I_{5,1}=T^2I_{5,1}$ is the union of a (full) subinterval of $I_{2,0}$ with a (left) subinterval of $I_{3,0}$, hence the ambiguous label). Thus for example $I_{7,1}$
is sent by $T$ onto $I_{7,0}$, then by another application of $T$ into $I_{1,0}$, hence $T_1=T^2$ on $I_{7,1}$. 
And we check that $T_1$ is indeed an AR9 interval exchange defined by  $(a_1,b_1,c_1)$ on the union of   $\Omega_1$, $\Omega'_1$, $\Omega"_1$; the order is still the first one.

\subsubsection{Induction step case I: $c_1=a_0-b_0-c_0.$} Then $a_1=b_0$, $b_1=c_0$.

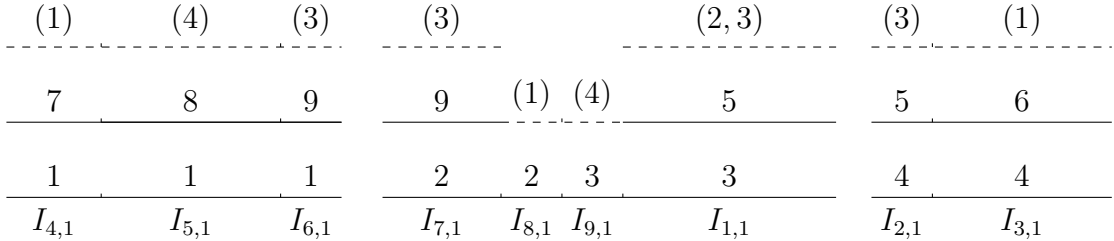
\begin{figure}[h] \label{Induction Case I}
\begin{center}
\begin{tikzpicture}[scale = 5]

\draw(0,.1)--(.891,.1);
\draw (.126,.1) node[above]{$1$};
\draw (.126,.1) node[below]{$I_{4,1}$};
\draw(.252,.1)--(.252,.11);
\draw (.49,.1) node[above]{$1$};
\draw (.49,.1) node[below]{$I_{5,1}$};
\draw(.729,.1)--(.729,.11);
\draw (.81,.1) node[above]{$1$};
\draw (.81,.1) node[below]{$I_{6,1}$};

\draw(1,.1)--(2.206,.1);
\draw (1.157,.1) node[above]{$2$};
\draw (1.157,.1) node[below]{$I_{7,1}$};
\draw(1.315,.1)--(1.315,.11);
\draw (1.396,.1) node[above]{$2$};
\draw (1.396,.1) node[below]{$I_{8,1}$};
\draw(1.477,.1)--(1.477,.11);
\draw (1.558,.1) node[above]{$3$};
\draw (1.558,.1) node[below]{$I_{9,1}$};
\draw(1.639,.1)--(1.639,.11);
\draw(1.922,.1) node[above]{$3$};
\draw(1.922,.1) node[below]{$I_{1,1}$};

\draw(2.3,.1)--(2.939,.1);
\draw (2.381,.1) node[above]{$4$};
\draw (2.381,.1) node[below]{$I_{2,1}$};
\draw(2.462,.1)--(2.462,.11);
\draw (2.7,.1) node[above]{$4$};
\draw (2.7,.1) node[below]{$I_{3,1}$};
\draw(2.939,.1) node[above]{};

\draw(1,.3)--(1.315,.3);
\draw(1.639,.3)--(2.206,.3);
\draw[dashed](1.315,.3)--(1.639,.3);
\draw (1.157,.3) node[above]{$9$};
\draw(1.315,.3) node[above]{};

\draw (1.396,.3) node[above]{$(1)$};
\draw(1.477,.3) node[above]{};
\draw(1.477,.3)--(1.477,.31);
\draw (1.558,.3) node[above]{$(4)$};
\draw(1.639,.3) node[above]{};

\draw(1.922,.3) node[above]{$5$};

\draw(2.3,.3)--(2.939,.3);
\draw (2.381,.3) node[above]{$5$};
\draw(2.462,.3) node[above]{};
\draw(2.462,.3)--(2.462,.31);
\draw (2.7,.3) node[above]{$6$};
\draw(2.939,.3) node[above]{};

\draw(0,.3)--(.891,.3);
\draw(.252,.3)--(.891,.3);
\draw (.126,.3) node[above]{$7$};
\draw(.252,.3) node[above]{};
\draw(.252,.3)--(.252,.31);
\draw (.49,.3) node[above]{$8$};
\draw(.729,.3) node[above]{};
\draw(.729,.3)--(.729,.31);
\draw (.81,.3) node[above]{$9$};
\draw(.891,.3) node[above]{};

\draw[dashed](1,.5)--(1.315,.5);
\draw[dashed](1.639,.5)--(2.206,.5);
\draw (1.157,.5) node[above]{$(3)$};
\draw(1.315,.5) node[above]{};

\draw (1.396,.5) node[above]{};
\draw(1.477,.5) node[above]{};

\draw (1.558,.5) node[above]{};
\draw(1.639,.5) node[above]{};

\draw(1.922,.5) node[above]{$(2,3)$};

\draw[dashed](2.3,.5)--(2.939,.5);
\draw (2.381,.5) node[above]{$(3)$};
\draw(2.462,.5) node[above]{};
\draw(2.462,.5)--(2.462,.51);
\draw (2.7,.5) node[above]{$(1)$};
\draw(2.939,.5) node[above]{};

\draw[dashed](0,.5)--(.891,.5);
\draw (.126,.5) node[above]{$(1)$};
\draw(.252,.5) node[above]{};
\draw(.252,.5)--(.252,.51);
\draw (.49,.5) node[above]{$(4)$};
\draw(.729,.5) node[above]{};
\draw(.729,.5)--(.729,.51);
\draw (.81,.5) node[above]{$(3)$};
\draw(.891,.5) node[above]{};

\end{tikzpicture}\\
\caption{Induction Case I}

\end{center}

\end{figure}

The length of each $I_{i,1}$ in Figure 5 is the same 
  as in case III. $T_1$ is an AR9 interval exchange defined by  $(a_1,b_1,c_1)$, in  the third order.

\subsubsection{Induction step case II: $b_1=a_0-b_0-c_0.$} Then $a_1=b_0$, $c_1=c_0$.

\begin{figure}[h] \label{Induction Case II}
\begin{center}
\begin{tikzpicture}[scale = 5]

\draw (0,.1) node[above]{};
\draw (1.35,.1) node[above]{};
\draw (2.05,.1) node[above]{};

\draw(0,.1)--(1.206,.1);
\draw (1.049,.1) node[above]{$1$};
\draw (1.049,.1) node[below]{$I_{7,1}$};
\draw(.891,.1) node[above]{};
\draw(.891,.1)--(.891,.11);
\draw (.81,.1) node[above]{$1$};
\draw (.81,.1) node[below]{$I_{8,1}$};
\draw(.729,.1) node[above]{};
\draw(.729,.1)--(.729,.11);
\draw (.648,.1) node[above]{$1$};
\draw (.648,.1) node[below]{$I_{9,1}$};
\draw(.567,.1) node[above]{};
\draw(.567,.1)--(.567,.11);
\draw(.284,.1) node[above]{$1$};
\draw(.284,.1) node[below]{$I_{1,1}$};

\draw(2.3,.1)--(2.939,.1);
\draw (2.858,.1) node[above]{$4$};
\draw (2.858,.1) node[below]{$I_{2,1}$};
\draw(2.777,.1) node[above]{};
\draw(2.777,.1)--(2.777,.11);
\draw (2.539,.1) node[above]{$4$};
\draw (2.539,.1) node[below]{$I_{3,1}$};
\draw(1.939,.1) node[above]{};

\draw(1.3,.1)--(2.191,.1);
\draw (2.065,.1) node[above]{$3$};
\draw (2.065,.1) node[below]{$I_{4,1}$};
\draw(1.939,.1) node[above]{};
\draw(1.939,.1)--(1.939,.11);
\draw (1.701,.1) node[above]{$3$};
\draw (1.701,.1) node[below]{$I_{5,1}$};
\draw(1.462,.1) node[above]{};
\draw(1.462,.1)--(1.462,.11);
\draw (1.381,.1) node[above]{$2$};
\draw (1.381,.1) node[below]{$I_{6,1}$};
\draw(2.891,.1) node[above]{};

\draw(0,.3)--(.891,.3);
\draw[dashed](.891,.3)--(1.206,.3);
\draw (1.049,.3) node[above]{$(1)$};
\draw(.891,.3) node[above]{};
\draw(.567,.3)--(.567,.31);
\draw (.81,.3) node[above]{$9$};
\draw(.729,.3) node[above]{};
\draw(.729,.3)--(.729,.31);
\draw (.648,.3) node[above]{$8$};
\draw(.567,.3) node[above]{};
\draw(.284,.3) node[above]{$7$};

\draw(2.3,.3)--(2.939,.3);
\draw (2.858,.3) node[above]{$6$};
\draw(2.777,.3) node[above]{};
\draw(2.777,.3)--(2.777,.31);
\draw (2.539,.3) node[above]{$5$};

\draw[dashed](1.3,.3)--(1.939,.3);
\draw(1.939,.3)--(2.191,.3);
\draw (2.065,.3) node[above]{$5$};
\draw(1.939,.3) node[above]{};
\draw (1.701,.3) node[above]{$(4)$};
\draw(1.462,.3) node[above]{};
\draw(1.939,.3)--(1.939,.31);
\draw (1.381,.3) node[above]{$(1)$};
\draw(2.891,.3) node[above]{};

\draw[dashed](0,.5)--(.891,.5);
\draw (.81,.5) node[above]{$(3)$};
\draw(.567,.5)--(.567,.51);
\draw (.648,.5) node[above]{$(4)$};
\draw(.729,.5)--(.729,.51);
\draw(.284,.5) node[above]{$(1)$};

\draw[dashed](1.939,.5)--(2.191,.5);
\draw (2.065,.5) node[above]{$(2,3)$};

\draw(2.777,.5)--(2.777,.51);
\draw[dashed](2.3,.5)--(2.939,.5);
\draw (2.858,.5) node[above]{$(1)$};
\draw (2.539,.5) node[above]{$(3)$};

\end{tikzpicture}\\
\caption{Induction Case II}

\end{center}

\end{figure}

The length of each $I_{i,1}$ in Figure 6 is the same as in case III.  $T_1$ is an AR9 interval exchange defined by  $(a_1,b_1,c_1)$, in  the reversed second order.\\

The same computations work if we start from an AR9  in the second order: we get the same pictures, in the second order in Case III, the first order in Case I, the reversed first order in Case II. When we start from the third  order, we get the same pictures, in
the third order in Case III, the reversed third order in Case II, and the second order in Case I. If we start form a reversed order, just reverse the orientation of the pictures.\\

We can now iterate the induction: starting  with $T_0=T$, we define $T_k$ as the induced map of $T_{k-1}$ on the set $\cup_{i=1}^4I_{i,k-1}$, which we denote by $J_{a,k-1}$. It defines  intervals $\Omega_k$, $\Omega'_k$, $\Omega"_k$, $I_{i,k}$, with   $\cup_{i=1}^9I_{i,k}=J_{a,k-1}$. The points $(a_k,b_k,c_k)$ in the Rauzy gasket have been defined in \cite{ar}, where the same induction process is described for  AR6 interval exchanges; they  constitute an algorithm of simultaneous approximation of $(a_0,b_0,c_0)$, which is called the {\em Arnoux-Rauzy algorithm}. Our induction on AR9 interval exchanges gives also an algorithm of simultaneous approximation of $(a_0+b_0,b_0+c_0,a_0+c_0)$ by the lengths of the intervals  $\Omega_k$, $\Omega'_k$, $\Omega"_k$: this turns out to be the 
{\em fully substractive algorithm} where the smallest of the numbers is substracted from the other two \cite{meno}.

\subsection{AR9 symbolic systems}\label{ccr}
\begin{definition}\label{dars} An AR9 symbolic system $(Y_9,S)$
 is the  {\em natural coding} of  an AR9 interval exchange $(X_9,T)$, that is its coding  by the partition into $I_{i,0}$, $1\leq i\leq 9$; we denote by $\psi$   the map  associating to each point $x \in X_9$ its trajectory in $Y_9$. \end{definition}
 
 \begin{remark}\label{rep} Because of the way we deal with the  endpoints, see Section \ref{nep} above, $\psi$ is injective but not surjective; we have $Y_9=\psi(X_9) \cup D_9$, where $D_9$ is a countable set made with the {\em improper trajectories} of the right endpoints of the intervals $I_{i,0}$ and their negative orbits: these are the limits, in the product topology of $\{1,...9\}^{\gN}$, in which $Y_9$ is closed, of trajectories of points approaching these endpoints from the left, and similarly for their pre-images.  \end{remark}

\begin{proposition}\label{ar9s} For each $(a_0,b_0,c_0)$ in $\Gamma$, the AR9 symbolic system associated to any AR9 interval exchange defined by $(a_0,b_0,c_0)$ is the  symbolic system on $\{1,... 9\}$ generated by the three substitutions 
 \begin{itemize}
 \item $\sigma'_I$:  $1\to 35$, $2\to 45$, $3\to 46$, $4\to 17$, $5\to 18$, $6\to 19$, $7\to 29$, $8\to 2$, $9\to 3$, 
  \item $\sigma'_{II}$:  $1\to 17$, $2\to 46$, $3\to 45$, $4\to 35$, $5\to 3$, $6\to 2$, $7\to 1$, $8\to 19$, $9\to 18$, 
   \item $\sigma'_{III}$:  $1\to 1$, $2\to 2$, $3\to 3$, $4\to 4$, $5\to 45$, $6\to 46$, $7\to 17$, $8\to 18$, $9\to 19$.
 \end{itemize}
and a directing sequence $r_n$,  $n\in \GN^{\star}$, $r_n \in\{I,II,III\}$, defined by $r_n=I$ if $a_n=a_{n-1}-b_{n-1}-c_{n-1}$, $r_n=II$ if $b_n=a_{n-1}-b_{n-1}-c_{n-1}$, $r_n=III$ if $c_n=a_{n-1}-b_{n-1}-c_{n-1}$; $r_n$ takes the value $I$ infinitely many times. 

Any system defined in this way is an AR9 symbolic system.\end{proposition}
{\bf Proof}\\
  We iterate the induction of Section \ref{sind}, and call  $i_k$, $1\leq i\leq 9$, the trajectory under $T$ of any point $x$ in $I_{i,k}$ between the time $0$ and the first return time of $x$ in $J_{a,k-1}$, coded by the partition into $I_{i,k}$, $1\leq i\leq 9$. The induction steps show that
 $i_k=\sigma'_{r_1}...\sigma'_{r_k}i$, and that, if we iterate the induction infinitely many times, the words $1_k$ to $9_k$, $k\geq 0$, generate the language of $T$. As $a_n>b_n>c_n>0$, $r_n=I$ infinitely often.
 
 It is actually proved in \cite{ar} that the construction of $r_n$ gives a one-to-one correspondence between the points of $\Gamma$ and the sequences $r_n$,  $n\in \GN^{\star}$, $r_n \in\{I,II,III\}$ where $r_n$ takes the value $I$ infinitely many times, which proves our last assertion.  \qed\\

Thus the AR9 symbolic system does not depend on the location or the order of the intervals $\Omega_0$, $\Omega'_0$, $\Omega"_0$.
 The common length of the words  $1_k$, $2_k$, $3_k$, $4_k$, is  $h_{a,k}$ defined in Section \ref{sar3}, $h_{b,k}$ is the common length of the words $5_k$, $6_k$, $7_k$, $h_{c,k}$ the common length of the words  $8_k$, $9_k$.\\

The multiplicative rules of Section \ref{pqm} above extend immediately to AR9 systems, in the following way
\begin{itemize}
 \item if the  $n+1$-th multiplicative rule is a rule  $I_m$,
 \begin{itemize} \item  $1_{m_{n+1}}=3_{m_n}4_{m_n}^{k_{n+1}-1}5_{m_n}$,
\item $2_{m_{n+1}}=4_{m_n}^{k_{n+1}}5_{m_n}$,
\item $3_{m_{n+1}}=4_{m_n}^{k_{n+1}}6_{m_n}$,
\item $4_{m_{n+1}}=1_{m_n}^{k_{n+1}}7_{m_n}$,
\item $5_{m_{n+1}}=1_{m_n}^{k_{n+1}}8_{m_n}$,
\item $6_{m_{n+1}}=1_{m_n}^{k_{n+1}}9_{m_n}$,
\item $7_{m_{n+1}}=2_{m_n}1_{m_n}^{k_{n+1}-1}9_{m_n}$,
 \item $8_{m_{n+1}}=2_{m_n}$,
\item $9_{m_{n+1}}=3_{m_n}$; \end{itemize}

 \item if the  $n+1$-th multiplicative rule is a rule  $II_m$,
\begin{itemize}\item $1_{m_{n+1}}=1_{m_n}^{k_{n+1}}7_{m_n}$,
\item $2_{m_{n+1}}=4_{m_n}^{k_{n+1}}6_{m_n}$,
\item $3_{m_{n+1}}=4_{m_n}^{k_{n+1}}5_{m_n}$,
 \item  $4_{m_{n+1}}=3_{m_n}4_{m_n}^{k_{n+1}-1}5_{m_n}$,
 \item $5_{m_{n+1}}=3_{m_n}$,
  \item $6_{m_{n+1}}=2_{m_n}$,
   \item $7_{m_{n+1}}=1_{m_n}$,
\item $8_{m_{n+1}}=1_{m_n}^{k_{n+1}}9_{m_n}$,
\item $9_{m_{n+1}}=1_{m_n}^{k_{n+1}}8_{m_n}$.
 \end{itemize}

 \end{itemize}

 \subsection{Relations between  Arnoux-Rauzy systems}\label{nno}
 Starting from a point $(a_0,b_0,c_0)$ in $\Gamma$,  we have defined two geometric systems, $(X_9,T)$ and $(X_6,T)$.  
 
  \begin{proposition}\label{ar63} An  AR9 interval exchange defined by $(a_0,b_0,c_0)$   is conjugate to an AR6 interval exchange defined by $(a_0,b_0,c_0)$ by a map which is continuous except on a finite number of points, and thus gives a measure-theoretic isomorphism for each invariant measure, and any AR6 interval exchange is conjugate to an AR9 in this way. \end{proposition}
{\bf Proof}\\
By gluing together the three intervals $\Omega_0$, $\Omega'_0$, $\Omega"_0$ we define a map $\phi'_6$ sending $X_9$ to a circle of length $2a_0+2b_0+2c_0$: for example, in the first order, we identify the right end of $\Omega_0$ with the left end of $\Omega'_0$,
the right end of $\Omega'_0$ with the left end of $\Omega'"_0$, the right end of $\Omega"_0$ with the left end of $\Omega_0$.  This conjugates $(X_9,T)$ to a system   $(X_6,T)$ which is  exactly the AR6 interval exchange  defined in Section \ref{sar6} above: its intervals of continuity are the $\phi'_6(J_{j,0})$, $j\in \{a-,a+,b-,b+,c-,c+\}$ where $J_{a-,0}=I_{1,0}\cup I_{2,0}$, $J_{a+,0}=I_{3,0}\cup I_{4,0}$, $J_{b-,0}=I_{5,0}$, $J_{b+,0}=I_{6,0}\cup I_{7,0}$, $J_{c-,0}=I_{8,0}$, $J_{c+,0}=I_{9,0}$. It is immediate that every AR6 interval exchange can be built in this way. \qed \\

 As in Proposition \ref{ar9s}, any point in $\Gamma$ defines a directing sequence $(r_n)$. Each directing sequence  defines two symbolic systems, $(Y_9,S)$ and $(Y_3,S)$.

\begin{proposition}\label{cod} The coding of an AR9 symbolic system defined by $(a_0,b_0,c_0)$, by the partition into three sets $J_{a,0}=I_{1,0}\cup I_{2,0}\cup I_{3,0}\cup I_{4,0}$, $J_{b,0}=I_{5,0}\cup I_{6,0}\cup I_{7,0}$, $J_{c,0}=I_{8,0}\cup I_{9,0}$, is the AR3 symbolic system defined by the directing sequence in Proposition \ref{ar9s}, and all AR3 symbolic systems can be built in this way. \end{proposition}
{\bf Proof}\\
 We define the letter-to-letter map $\phi$ by $\phi (1)=\phi (2)=\phi (3)=\phi (4)=a$, $\phi (5) =\phi (6)=\phi (7)=b$, $\phi (8)=\phi (9)=c$. If we build the words $A_k$, $B_k$, $C_k$  in Definition \ref{dar3} with a directing sequence $(r_n)$  and the words $1_k$ to $9_k$ in the proof of Proposition \ref{ar9s}, we get inductively that for all $k$, 
 $\phi (1_k)=\phi (2_k)=\phi (3_k)=\phi(4_k)=A_k$, $\phi (5_k) =\phi (6_k)=\phi (7_k)=B_k$, $\phi (8_k)=\phi (9_k)=C_k$. By the induction steps of Section \ref{sind}, $A_k$, resp. $B_k$, resp. $C_k$ are the trajectory, under the coding by the partition into three sets,  of any point $x$ in $I_{i,k}$, $1\leq i\leq 4$,  resp. $5\leq i\leq 7$, resp. $8\leq i\leq 9$, between the time $0$ and the first return time of $x$ in $J_{a,k-1}$, and the words $A_k$, $B_k$, $C_k$, $k\geq 0$, generate the language of $T$.
 
 The last assertion comes again from the fact, known from \cite{ar}, that each directing sequence defines a point in $\Gamma$. \qed \\

\begin{figure}[h] \label{x}
\begin{center}
	\begin{tikzpicture}[every text node part/.style={align=center}]
		\node (T) at (0,10) {$(X_6,T)$};
		\node (S) at (8,14) {$(Y_3,S)$};
		\node (M) at (0,14) {$(Y_6,S)$};
			\node (B) at (-8,10) {$(X_9,T)$};
		\node (TL) at (-8,14) {$(Y_9,S)$};
		
	\draw[->] (T) edge node[auto] {$\psi'$} (M);		
		\draw[->] (B) edge node[auto] {$\phi'_6$} (T);	
		\draw[->,dashed] (TL) edge[bend left=30] node[auto] {$\phi$} (S);	 
		\draw[->] (B) edge node[auto] {$\psi$} (TL);	 
		\draw[->] (TL) edge node[auto] {$\phi_6$} (M);
			\draw[->,dashed] (M) edge node[auto] {$\phi_3$} (S);

	\end{tikzpicture}\\
	\caption{The five AR systems}

\end{center}

\end{figure}
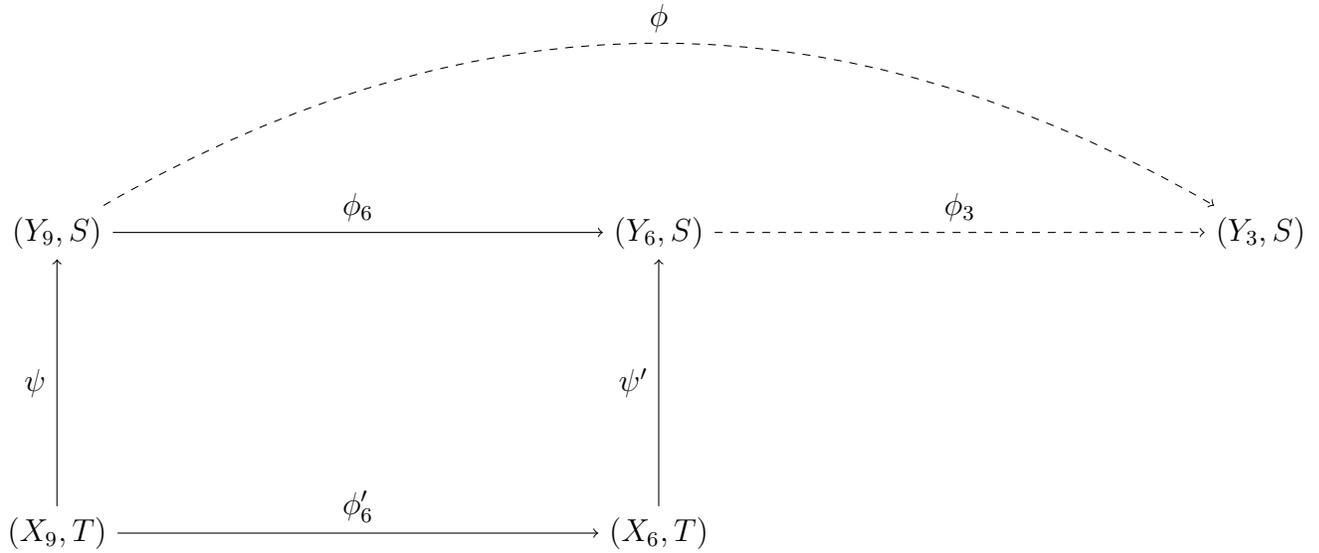

\begin{corollary} An  AR9 symbolic system  has an AR3 symbolic system defined by the same directing sequence as a factor, and all AR3 symbolic systems can be built in this way. \end{corollary}
 {\bf Proof}\\
 These are two codings of the same AR9 interval exchange, and the partition into nine intervals is finer than the partition into three sets. \qed\\

The map associating to a point in $(X_9,S)$ its coding in $(Y_3,S)$ is just $\phi\psi$, where $\psi$ is defined in Definition \ref{dars} and $\phi$ in the proof of Proposition \ref{cod}. As in Remark \ref{rep}, we have $\phi\psi (X_9)=Y_3 \setminus D'_3$ for the countable set $D'_3$ made with improper trajectories; note that $D'_3\subset D_3$ where $D_3=\phi(D_9)$. $\phi\psi$ conjugates the map $T$ on $X_9$ with the shift $S$ on $X_3$: to use the vocabulary of \cite{ar}, $\phi\psi$ is called a {\em semi-conjugacy}; as is pointed out in the introduction above, this does not give a one-to-one correspondence between points. Similarly, $\phi$ conjugates the shifts on $Y_9$ on $Y_3$ and $\phi (Y_9)=Y_3$; it is also a semi-conjugacy, and not injective, see Proposition \ref{E123}  below.\\

We can also define the {\em AR6 symbolic system}  $(Y_6,S)$ on $\{a-,a+,b-,b+,c-,c+\}$, as the natural coding $\psi_6$, of $(X_6,T)$, by its  six intervals of continuity; we have $Y_6=\psi_6(X_6)\cup D_6$ for a countable set $D_6$. We can write $\phi=\phi_3\circ\phi_6$, with $\phi_6 (1)=\phi_6 (2)=a-$, $\phi_6 (3)=\phi_6 (4)=a+$, $\phi_6 (5) =b-$, $\phi_6 (6)=\phi_6 (7)=b+$, $\phi_6 (8)=c-$, $\phi_6 (9)=c+$, and $\phi_{3}(j-)=\phi_{3}(j+)=j$ for $j=a,b,c$. 

In the same way as Proposition \ref{cod}, we could reprove the main result of \cite{ar}: 
the coding of an AR6 interval exchange defined by $(a_0,b_0,c_0)$, by the partition into three sets $\phi'_6(J_{a-,0}\cup J_{a+,0})$, $\phi'_6(J_{b-,0}\cup J_{b+,0})$, $\phi'_6(J_{c-,0}\cup J_{c+,0})$, is the AR3 symbolic system defined by the directing sequence of Proposition \ref{ar9s}, and all AR3 symbolic systems can be built in this way. Thus $(Y_6,S)$ appears as an intermediate coding between the AR3 and AR9 symbolic systems; because of Proposition \ref{ar63}, $\phi_6$, applied letter to letter, is invertible except on a countable set (included in $\phi_6(D_9)$), and conjugates $(Y_9,S)$ and $(Y_6,S)$, which are thus measure-theoretically isomorphic for each invariant measure.

In Figure 7, the four systems linked by full edges are, for all our purposes, the same system; the nature of the dashed edges will be investigated in the remainder of this paper.

As was already mentioned, we do not know any way to build the trajectories in $Y_6$ as in Definition \ref{dar3} or Proposition \ref{ar9s}; but they can be deduced from the  trajectories in $Y_9$ by applying $\phi_6$ letter to letter, and that was the main objective of the theory of  AR9 systems; however, in general it will be easier to work directly on AR9 systems and then derive the properties of AR6 systems. \\

At this stage, it may be useful to recall the various notations we use, for which we had to make choices because of the number of systems we have defined and some long pre-existing notations:
$a$, $b$, $c$ are always the three symbols on which AR3 systems are built. But $a_k$, $b_k$, $c_k$, for any $k$, are real numbers, describing  lengths of intervals. $A_k$, $B_k$, $C_k$ are the words used to build AR3 systems, of lengths (i.e. number of letters) $h_{a,k}$, $h_{b,k}$, $h_{c,k}$. $1$ to $9$ are the symbols  on which AR9 symbolic systems are built, $1_k$ to  $9_k$ are the words used to build them, their lengths are among $h_{a,k}$, $h_{b,k}$, $h_{c,k}$. Interval lengths for AR9 systems, when needed,  are  defined in terms of $a_k$, $b_k$, $c_k$. Roman numerals are used to number substitutions and rules to build words. 

\section{Dynamical properties}
\subsection{Minimality}
By using the condition that $r_n=I$ for infinitely many $n$, the minimality of AR3 symbolic systems and AR6 interval exchanges is shown in \cite{ar}. The minimality of AR6 symbolic systems follows, as the minimality of an interval exchange is equivalent to the minimality of its natural coding, small intervals corresponding to small cylinders. 

\begin{proposition} Any AR9 system  is minimal. \end{proposition}
{\bf Proof}\\
We show it for the symbolic systems, the minimality of the interval exchanges follows from the remark just above. We want to show that in the language  of $(Y_9,S)$ any word $w$ occurs in any long enough word. It is enough to show that for all $n$ and $1\leq i\leq 9$ there exists $N$ such that $i_n$ occurs in every $j_N$, $1\leq j\leq 9$.

For example we take $i=1$. Through  $\sigma'_{III}$ $i_n$ occurs in $i_{n+1}$ for all $i$, as we are after sufficient conditions we can ignore these  rules. We start from $1_n$; it occurs  in $1_{n+1}$ through any number of  $\sigma'_{III}$, so we wait until the first  $\sigma'_I$ (we know it exists), in which $1_{p_1}$ occurs in $4_{p_1+1}$, $5_{p_1+1}$, $6_{p_1+1}$. \\

We follow these three words until just before the next $\sigma'_{I}$: if there has been no  $\sigma'_{II}$, we have to track $4_{p_2}$, $5_{p_2}$, $6_{p_2}$; if there has been one $\sigma'_{II}$, the words into which at least one of $4_{p_1+1}$, $5_{p_1+1}$, $6_{p_1+1}$ occur are $2_{p_2}$, $3_{p_2}$, $4_{p_2}$; if there have been two $\sigma'_{II}$ or more,  these words are $2_{p_2}$, $3_{p_2}$, $4_{p_2}$, $5_{p_2}$, $6_{p_2}$. So in the worst case we have to track either $2_{p_2}$, $3_{p_2}$, $4_{p_2}$ or $4_{p_2}$, $5_{p_2}$, $6_{p_2}$. After the $\sigma'_{I}$, these occur either in $1_{p_2+1}$, $2_{p_2+1}$, $3_{p_2+1}$ or in a larger set of words. 

Again we follow these three words until just before the next  $\sigma'_I$: if there have been no $\sigma'_{II}$, we have to track $1_{p_3}$, $2_{p_3}$, $3_{p_3}$; if there has been one $\sigma'_{II}$, the words to track are $1_{p_3}$, $4_{p_3}$, $5_{p_3}$, $6_{p_3}$, $7_{p_3}$, $8_{p_3}$, $9_{p_3}$; if there have been two $\sigma'_{II}$,   these words are $1_{p_3}$, $2_{p_3}$, $3_{p_3}$, $4_{p_3}$, $7_{p_3}$, $8_{p_3}$, $9_{p_3}$; if there have been at least three $\sigma'_{II}$, we have already won ($i_n$ occurs in all the $j_{p_3}$).

$1_{p_3}$, $4_{p_3}$, $5_{p_3}$, $6_{p_3}$, $7_{p_3}$, $8_{p_3}$, $9_{p_3}$ after $\sigma'_{I}$ give $1_{p_4}$, $2_{p_4}$, $3_{p_4}$, $4_{p_4}$, $5_{p_4}$, $6_{p_4}$, $7_{p_4}$ which are conserved by any number of $\sigma'_{II}$, and give every word after the next  $\sigma'_{I}$.\\
$1_{p_3}$, $2_{p_3}$, $3_{p_3}$, $4_{p_3}$, $7_{p_3}$, $8_{p_3}$, $9_{p_3}$ give everything after  $\sigma'_{I}$.\\
$1_{p_3}$, $2_{p_3}$, $3_{p_3}$ after $\sigma'_I$ give $1_{p_4}$, $4_{p_4}$, $5_{p_4}$, $6_{p_4}$, $7_{p_4}$, $8_{p_4}$, $9_{p_4}$ (with which we win after another $\sigma'_{I}$, as just above), after one $\sigma'_{II}$ $1_{p_4}$, $2_{p_4}$, $3_{p_4}$, $4_{p_4}$, $7_{p_4}$, $8_{p_4}$, $9_{p_4}$  which will give everything after $\sigma'_{I}$,  after two $\sigma'_{II}$ everything.\\

Similar (shorter, as we can use what we already proved about $1_n$ and successive others) chasing arguments take care of the other $i_n$.\qed\\

\subsection{Rokhlin towers}

\begin{definition}  In a system $(X',U)$, a  {\em Rokhlin tower}  is a collection 
    of disjoint measurable sets $F$, $UF$, \ldots, $U^{h-1}F$ ($U^jF$ is called {\em level} $j$ of the tower, $F$ is called the {\em base}, $h$ the {\em height} of the tower).  A {\em slice} of $\tau$ is a union of whole levels $U^{p_1}F$ ... $U^{p_l}F$, and a {\em column} of $\tau$ is  a union of all sublevels $G$, ... $U^{h-1}G$ for a subset $G$ of $F$. We shall usually write ``the tower $\tau$" as a shortened form of ``the tower  for which the union of the levels is the set $\tau$".\end{definition}

\begin{proposition}\label{tto}
In an AR9 interval exchange $(X_9,T)$, there are nine sequences of towers $\tau_{i,k}$, respectively of base $I_{i,k}$, and height equal to the length of the word $i_k$, $1\leq i \leq 9$, $k\geq 0$: the union of all the levels for fixed $k$ is $X_9$, and every point $x$ in $X_9$ is determined by the sequence $\iota(x,k)$, $\eta(x,k)$ such that $x$ is in $T^{\eta(x,k)}I_{\iota(x,k),k}$, $k\geq 0$. This remains true if we restrict $k$ to a subsequence, for example the $m_n$. All levels of these  towers are intervals. \end{proposition}
{\bf Proof}\\
From the induction steps in Section \ref{sind}, we deduce that the $\tau_{i,k}$ are indeed Rokhlin towers, whose  union of levels for fixed $k$ is indeed $X_9$; all these levels  are intervals and  their lengths are smaller than $a_k$, which tends to zero when $k$ goes to infinity, hence the result. \qed\\

Figures $4$, $5$, $6$ going from stage $0$ to stage $1$
show how the towers at order $1$ are made from the towers at order $0$ by {\em cutting and stacking}. This cutting and stacking is done in the same way from stage $k$ to stage $k+1$; it is dictated by the induction as above, and can be read on the rules giving the words $1_{k+1}$ to $9_{k+1}$ as concatenations of the words $1_k$ to $9_k$, which are deduced  from the substitutions $\sigma'_I$ to $\sigma'_{III}$: for example, when $r_{k+1}=I$, $\sigma'_I$ is applied, and we deduce from $1\to 35$ that $1_{k+1}=3_k5_k$, and the tower $\tau_{1,k+1}$ is made by a column of $\tau_{5,k}$  stacked above a column of $\tau_{3,k}$. \\

\begin{corollary} In $(Y_9,S)$, the $\tau'_{i,k}=\psi(\tau_{i,k})$, $i=1,...9$, form nine sequences of Rokhlin towers. If $D_9$ is the countable set defined in Remark \ref{rep}, every point $y$ in $Y_9\setminus D_9$ is determined by the sequences $\iota(y,k)$, $\eta(y,k)$ such that $y$ is in $S^{\eta(x,k)}\psi(I_{\iota(x,k),k})$, $k\geq 0$.

In $(X_9,T)$, there exist three sequences of Rokhlin towers $\tau_{a,k}$, $\tau_{b,k}$, $\tau_{c,k}$, respectively of bases $J_{a,k}$, $J_{b,k}$, $J_{c,k}$, and heights  equal to $h_{a,k}$, $h_{b,k}$, $h_{c,k}$, $k\geq 0$, where 
 $J_{a,k}=I_{1,k}\cup I_{2,k}\cup I_{3,k}\cup I_{4,k}$, 
 $J_{b,k}=I_{5,k}\cup I_{6,k}\cup I_{7,k}$,  $J_{c,k}=I_{8,k}\cup I_{9,k}$. The union of all their levels for fixed $k$ is $X_9$. 
 
 In the AR3 system $(Y_3,S)$, the $\tau'_{j,k}=\phi\psi(\tau_{j,k})$, $j=a,b,c$, form three sequences of Rokhlin towers; if  $D_3=\phi(D_9)$, every point $x$ in $Y_3\setminus D_3$ is determined by the sequences $\iota'(y,k)$, $\eta(y,k)$ such that $y$ is in $S^{\eta(y,k)}\phi\psi(J_{\iota'(y,k),k})$, $k\geq 0$.\end{corollary}
 {\bf Proof}\\
 The first assertion comes from Proposition \ref{tto} translated by $\psi$ to the symbolic system, the second one from the definition of the $J_{j,k}$ and the values of the heights, the third one from the first one and the fact that 
 for all $k$  $\phi$ sends  $\psi(I_{i,k})$ to $\psi(J_{a,k})$ if $i=1,2,3,4$, $\psi(J_{b,k})$ if $i=5,6,7$, $\psi(J_{c,k})$ if$i=8,9$. and similarly for the other levels. \qed\\

\begin{remark}\label{r2} We can also build directly  (slightly) enlarged versions of the various towers $\tau'$ in the symbolic systems:  this is done  in \cite{cfme} for the  $\tau'_{j,k}$, $j=a,b,c$,  by induction on  cylinders which are the closure of $\phi\psi(J_{a,k})$ in the topology of the symbolic systems, and can be done in the same way for the  $\tau'_{i,k}$, $i=1,...9$,  by induction on  unions of cylinders which are the closure of  $\psi(J_{a,k})$. These enlarged towers are closed  and include also improper trajectories; but we do not need that for our results, for which countable sets can be neglected, and in any case points of $D_3$ must be taken into account, see Remark \ref{rd} below.\end{remark}

 The towers  $\tau'_{i,k}$, $i=1,...9$,  can be built by cutting and stacking  with the same rules as the  $\tau_{i,k}$. The $\tau_{j,k}$ or $\tau'_{j,k}$, $j=a,b,c$,  can be built by cutting and stacking, using the concatenation rules generating the words $A_k$ to $C_k$, deduced from the substitutions $\sigma_I$ to $\sigma_{III}$; we shall also use the multiplicative rules to build more quickly these towers at multiplicative times, as is shown in Figures 9 and 10 below.\\
 
 \begin{lemma}\label{toursint} For every $k$, the sets $T^jI_{2,k}$ and $T^jI_{3,k}$, $0\leq j \leq h_{a,k}-1$,  resp. $T^jI_{5,k}$ and $T^ji_{6,k}$, $0\leq j \leq h_{b,k}-1$,  resp. $T^jI_{8,k}$ and $T^ji_{9,k}$, $0\leq j \leq h_{c,k}-1$, are  adjacent  intervals. \end{lemma}
{\bf Proof}\\
We make the induction hypothesis that our result is true at order $k$ and that $T^jI_{2,k}$,  $T^jI_{5,k}$, r $T^jI_{8,k}$ are the leftmost of the respective two adjacent intervals when $T_k$ is not in a reversed order, the rightmost if $T_k$ is in a reversed order.

This is true for $k=0$, whatever the order. 
The induction step from $k$ to $k+1$ describes also the way the towers at order $k+1$ are built from the towers at order $k$. \\

Take for example Case I when $T_k$ is  not in a reversed order: the new tower $8$ is  made by taking a right subinterval of the base $I_{2,k}$ of the old tower $2$, and keeping the corresponding part of all the levels of the old tower $2$; the new tower $9$ is made by taking a left subinterval of the base $I_{3,k}$ of the old tower $3$, and keeping the corresponding part of all the levels of the old tower $3$. Thus all corresponding levels of the new towers $8$ and $9$ are adjacent as those of the old towers $2$ and $3$ were, and the levels of the new tower $8$ are to the left of those of the new tower $9$. \\
The new tower $2$ is  made by taking a left subinterval of the base $I_{4,k}$ of the old tower $4$, and keeping the corresponding part of all the levels of the old tower $4$, until the top; above that we stack  a right subinterval of  $I_{5,k}$,  and the corresponding part of all the levels of the old tower $5$. The new tower $3$ is  made by taking a right subinterval of $I_{4,k}$, and keeping the corresponding part of all the levels of the old tower $4$, until the top; above that we stack $I_{6,k}$,  and  all the levels of the old tower $6$. Thus all corresponding levels of the new towers $2$ and $3$ are adjacent as those of the old towers $5$ and $6$ were, while the levels of the old tower $4$ are intervals, and the levels of the new tower $2$ are to the left of those of the new tower $3$. \\
The new tower $6$ is  made by taking a right subinterval of  $I_{1,k}$, and keeping the corresponding part of all the levels of the old tower $1$, until the top; above that we stack  a left subinterval of  $I_{9,k}$,  and the corresponding part of all the levels of the old tower $9$. The new tower $5$ is  made by taking a subinterval of $I_{1,k}$ just left of the previous one, and keeping the corresponding part of all the levels of the old tower $1$, until the top; above that we stack  the subinterval  $I_{8,k}$,  and all the levels of the old tower $8$. Thus all corresponding levels of the new towers $5$ and $6$ are adjacent as those of the old towers $8$ and $9$ were, while the levels of the old tower $1$ are intervals, and the levels of the new tower $5$ are to the left of those of the new tower $6$. \\

The other cases are similar.\qed

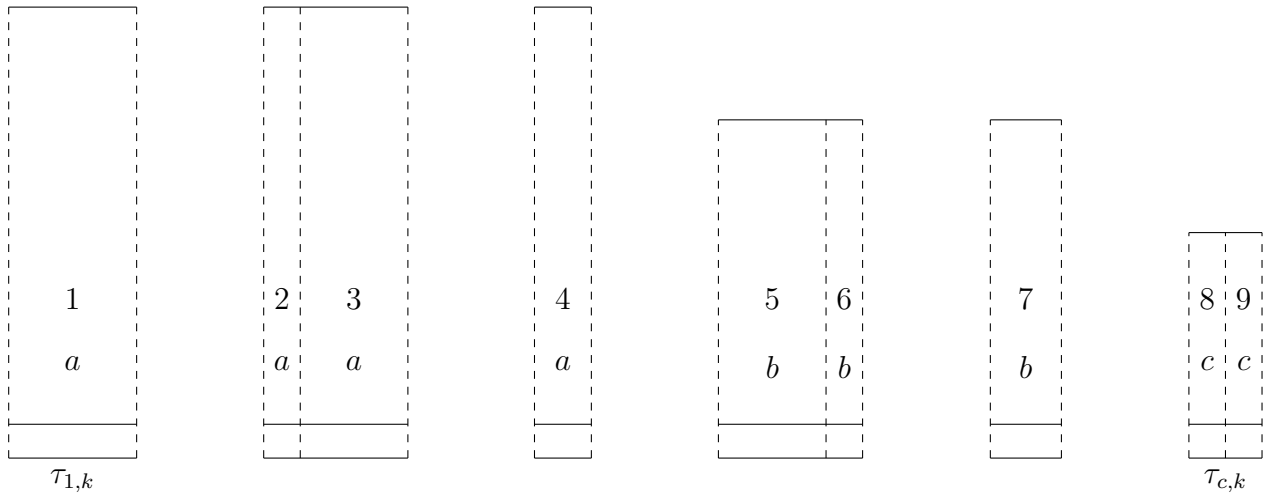
\begin{figure}[h] \label{Towers}
\begin{center}
\begin{tikzpicture}[scale = 3]

\draw(0,.1)--(.567,.1);
\draw(1.130,.1)--(1.769,.1);
\draw(2.330,.1)--(2.582,.1);
\draw(3.145,.1)--(3.784,.1);
\draw(4.35,.1)--(4.665,.1);
\draw(5.23,.1)--(5.554,.1);

\draw(0,.25)--(.567,.25);
\draw(1.130,.25)--(1.769,.25);
\draw(2.330,.25)--(2.582,.25);
\draw(3.145,.25)--(3.784,.25);
\draw(4.35,.25)--(4.665,.25);
\draw(5.23,.25)--(5.554,.25);

\draw(0,2.1)--(.567,2.1);
\draw(1.130,2.1)--(1.769,2.1);
\draw(2.330,2.1)--(2.582,2.1);
\draw(3.145,1.6)--(3.784,1.6);
\draw(4.35,1.6)--(4.665,1.6);
\draw(5.23,1.1)--(5.554,1.1);

\draw[dashed](3.145,.1)--(3.145,1.6);
\draw[dashed](3.622,.1)--(3.622,1.6);
\draw[dashed](3.784,.1)--(3.784,1.6);
\draw[dashed](4.35,.1)--(4.35,1.6);
\draw[dashed](4.665,.1)--(4.665,1.6);

\draw[dashed](5.23,.1)--(5.23,1.1);
\draw[dashed](5.392,.1)--(5.392,1.1);
\draw[dashed](5.554,.1)--(5.554,1.1);

\draw(.283,.6) node[below]{$a$};
\draw(.283,.9) node[below]{$1$};
\draw(1.211,.6) node[below]{$a$};
\draw(1.211,.9) node[below]{$2$};
\draw(1.53,.6) node[below]{$a$};
\draw(1.53,.9) node[below]{$3$};
\draw(2.456,.6) node[below]{$a$};
\draw(2.456,.9) node[below]{$4$};
\draw(3.384,.6) node[below]{$b$};
\draw(3.384,.9) node[below]{$5$};
\draw(3.703,.6) node[below]{$b$};
\draw(3.703,.9) node[below]{$6$};
\draw(4.507,.6) node[below]{$b$};
\draw(4.507,.9) node[below]{$7$};
\draw(5.311,.6) node[below]{$c$};
\draw(5.311,.9) node[below]{$8$};
\draw(5.473,.6) node[below]{$c$};
\draw(5.473,.9) node[below]{$9$};
\draw(.283,.1) node[below]{$\tau_{1,k}$};
\draw(5.392,.1) node[below]{$\tau_{c,k}$};

\draw[dashed](0,.1)--(0,2.1);
\draw[dashed](.567,.1)--(.567,2.1);
\draw[dashed](1.130,.1)--(1.130,2.1);
\draw[dashed](1.292,.1)--(1.292,2.1);
\draw[dashed](1.769,.1)--(1.769,2.1);
\draw[dashed](2.33,.1)--(2.33,2.1);
\draw[dashed](2.582,.1)--(2.582,2.1);

\end{tikzpicture}\\
\caption{Rokhlin towers in $X_9$}

\end{center}

\end{figure}

An immediate consequence is best seen on Figure 8:

\begin{corollary}\label{tourscor} Each level of the towers $\tau_{c,k}$ is an interval, each level of the towers $\tau_{b,k}$ is a  union of at most  two intervals, each level of the towers $\tau_{a,k}$ is a union of at most three intervals. \end{corollary}

 Note that  the $J_{j,k}$ and their images are not  intervals for $j=a,b$, except maybe for the first values of $k$, with some choices of $\Omega_0$, $\Omega'_0$, $\Omega"_0$, but even in that case, for example if they are in the first order, $J_{b,0}$ is not an interval. Similarly, except maybe for the first values of $k$, the levels of the towers $\tau_{b,k}$ are not intervals, the levels of the towers $\tau_{a,k}$ are not unions of less than three intervals.

\subsection{Isomorphism}\label{simp}

 \begin{definition} For $i=1,2,3$, let $E_i \subset Y_3$ be the set of points which have $i$ pre-images under $\phi$. \end{definition}

\begin{proposition}\label{E123} $Y_3\setminus D_3\subset E_1\cup E_2\cup E_3$. $E_3$ is countable. If $\mu(E_1)<1$, then for any invariant probability $\mu'$ the system $(Y_9,S, \mu')$  is a two-point extension of $(Y_3,S,\mu)$ \end{proposition}
{\bf Proof}\\ 
The knowledge of a point $y$ in $Y_3$ determines the sequences $\iota'(y,k)$ in $\{a,b,c\}$ and $0\leq \eta
(y,k)\leq h_{\iota'(y,k),k}-1$ such that $y$ is in $S^{\eta(y,k)}\phi\psi J_{\iota'(y,k),k}$ for all $k$. Except if  $y$ is in $D_3$, there exist points in $x\in X_9$ such that   $\phi\psi(x)=y$, and the  pre-images of $y$ by $\phi$ are the points $\psi(x)$; because of the way $\phi\psi$ acts on the towers, all these $x$ must be in $T^{\eta(y,k)}I_{\iota(x,k),k}$ were $\iota'(y,k)=\phi (\iota(x,k))$. $y$ being fixed, for a given $k$, all possible $x$ are in at most three of the intervals of Figure 8 above: if $\iota'(y,k)=a$, all possible $x$ are either in $T^{\eta(y,k)}I_{1,k}$, or in $T^{\eta(y,k)}I_{2,k}\cup T^{\eta(y,k)}I_{3,k}$, or in $T^{\eta(y,k)}I_{4,k}$, and similarly there are only two possible intervals if  $\iota'(y,k)=b$, and one if  $\iota'(y,k)=c$. If there exist more than three such points $x$, two of them must be infinitely often in the same interval, thus must be the same as the intersection of infinitely many of these intervals defines at most one point. Thus we get our first assertion. 

By the same reasoning, if  $y\in Y_3\setminus D_3$ is in $\tau'_{c,k}\cup \tau'_{b,k}$ for infinitely many $k$, then $y$ is in $E_1\cup E_2$. Thus if $y$ is in $E_3\setminus D_3$, $y$ is in $\tau'_{a,k}$ for all $k\geq k_0$. By the rules of construction by cutting and stacking, this implies that for all $k\geq k_0$ $\eta(y,k)$ takes the same value $\eta_0$, thus any pre-image of $y$ by $\phi\psi$ is in $\cap_{k\geq k_0} T^{\eta_0}J_{a,k}$. For $\eta_0=0$, it is shown in \cite{ar} that this intersection consists indeed of three distinct points, whose images by $\psi$ are not in $D_6$ and which have the same images by $\phi\psi$, thus $E_3\setminus D_3$ consists of the union of the positive orbits of these three points, which proves our second assertion.

Thus $\mu(E_1\cup E_2)=1$, and if $\mu(E_1)<1$
 the number of pre-images by $\phi$ is two on a set of positive measure, thus almost everywhere by ergodicity, and this is our third assertion.\qed\\

\begin{lemma}\label{tourc} Let $y$ be in $Y_3\setminus D_3$. If $y$ is in $\tau'_{c,k}$  for infinitely many $k$, then $y$ is in $E_1$. \end{lemma}
{\bf Proof}\\
 Under the hypothesis, as in the proof of Proposition \ref{E123}, for infinitely many $k$ all the pre-images of $y$ by $\phi\psi$ are in an interval, of length $c_k$, thus the intersection of infinitely many of these intervals defines at most one point. \qed\\

\begin{remark}\label{rd} If  we  enlarge the  towers to cover all $Y_3$ as in \cite{cfme} and Remark \ref{r2} above, the generalization of Lemma \ref{tourc} does not hold for $y\in D_3$: indeed,  the point $x_0$ separating $I_{8,0}$ from $I_{9,0}$ defines one trajectory in $\psi(X_9)$ and one improper trajectory (as in Remark \ref{rep}), and both these trajectories have the same image $y_0$ by $\phi$, though we can check that, for example in the Tribonacci case, $y_0$ is in the enlarged $\tau'_{c,k}$ for infinitely many $k$. However, it is true that every point in $Y_3$ has at most three pre-images by $\phi$, as the only candidates to have more are the points which are  in the enlarged $\tau'_{a,k}$ for all $k\geq k_0$, and their pre-images do not give rise to improper trajectories. \end{remark}

At this stage, one can ask whether the condition to be in  $\tau'_{c,k}$ for infinitely many $k$ is necessary for $y$ to be in $E_1$. Hopefully, a necessary and sufficient condition will be given in a further paper, but the following lemma gives already a negative answer for many systems including Tribonacci.

\begin{lemma}\label{tourab} Suppose that, 
 \begin{itemize} 
\item{(i)} either  for an infinite sequence $s_j$, the  $s_j+2$-th  multiplicative rule is $I_m$ with $k_{s_j+2}=1$, \item{(ii)} or for an infinite sequence $s_j$ the  $s_j+2$-th  multiplicative rules is $I_m$ and the $s_j+1$-th multiplicative rule is $II_m$ with $k_{s_j+1}=1$.
\end{itemize} Let $y$ be in $Y_3\setminus D_3$. If we are in  case $(i)$ and for infinitely many $j$ $y$ is in $\tau'_{b,m_{s_j+1}}\cap\tau'_{b,m_{s_j+3}}$, 
or if we are in  case $(ii)$ and for infinitely many $j$ $y$ is in $\tau'_{b,m_{s_j}}\cap\tau'_{b,m_{s_j+3}}$, then $y$ is in $E_1$. \end{lemma}
{\bf Proof}\\
A pre-image $x$ of  $y$ by $\phi$ s in $\tau'_{5,m_{s_j+3}}$, $\tau'_{6,m_{s_j+3}}$, or $\tau'_{7,m_{s_j+3}}$.

Going from $m_{s_j+2}$ to $m_{s_j+3}$, we have a number (possibly zero) of  $\sigma'_{III}$ followed by a  $\sigma'_{I}$ or $\sigma'_{II}$.
\begin{itemize}\item  Suppose this last substitution is $\sigma'_{II}$: the construction of the towers by  $\sigma'_{II}$ implies that $x$ is in $\tau'_{1,m_{s_j+3}-1}$, $\tau'_{2,m_{s_j+3}-1}$,  or $\tau'_{3,m_{s_j+3}-1}$; then either the absence of $\sigma'_{III}$ or the construction of the towers by $\sigma'_{III}$ imply that $x$ is in $\tau'_{1,p}$, $\tau'_{2,p}$ or $\tau'_{3,p}$ at all  stages $m_{s_j+2}\leq p \leq ,m_{s_j+3}-1$. 
\item Suppose now this substitution is $\sigma'_{I}$: the construction of the towers by  $\sigma'_{I}$ implies that $x$ is in $\tau'_{1,m_{s_j+3}-1}$, $\tau'_{2,m_{s_j+3}-1}$, $\tau'_{8,m_{s_j+3}-1}$ or $\tau'_{9,m_{s_j+3}-1}$. In the last two cases, $x$ is in $\tau'_{c,m_{s_j+3}-1}$ and if this happens infinitely often we conclude by Lemma \ref{tourc} that $y$ is in $E_1$. Otherwise, either the absence of $\sigma'_{III}$ or the construction of the towers by $\sigma'_{III}$ imply that $x$ is in $\tau'_{1,p}$ or $\tau'_{2,p}$ at all  stages $m_{s_j+2}\leq p\leq m_{s_j+3}-1$. 
\end{itemize}
Thus in both remaining cases $x$ is in $\tau'_{1,m_{s_j+2}}$, $\tau'_{2,m_{s_j+2}}$, or $\tau'_{3,m_{s_j+2}}$.

 Going from $m_{s_j+1}$ to $m_{s_j+2}$, we have a number of $\sigma'_{III}$ followed by a $\sigma'_{I}$; the construction of the towers by  $\sigma'_{I}$ implies that $x$ is in $\tau'_{3,m_{s_j+2}-1}$, $\tau'_{4,m_{s_j+2}-1}$, $\tau'_{5,m_{s_j+2}-1}$, or $\tau'_{6,m_{s_j+2}-1}$.  We are in the last two cases whenever $x$ is in $\tau'_{b,m_{s_j+2}-1}$, and then the  knowledge of its level in that tower puts $x$ in a single level of $\tau'_{5,m_{s_j+2}-1}\cup \tau'_{6,m_{s_j+2}-1}$, which puts the possible pre-images of $y$ by $\phi\psi$ in a small interval by Lemma \ref{toursint}; if this happens infinitely often we conclude as in  Lemma \ref{tourc} that $y$ is in $E_1$. Otherwise, either the absence of $\sigma'_{III}$ or the construction of the towers by $\sigma'_{III}$ imply that $x$ is in $\tau'_{3,p}$ or $\tau'_{4,p}$ at all  stages $m_{s_j+1}\leq p\leq ,m_{s_j+2}-1$: this is excluded by the hypotheses in case $(i)$, thus our result in proved in that case.

Finally, in case $(ii)$, going from $m_{s_j}$ to $m_{s_j+1}$ by  a single $\sigma'_{II}$ and knowing $y$ is in $\tau'_{b,m_{s_j}}$, we get that $x$ must be in   $\tau'_{5,m_{s_j}}$, and the  knowledge of its level in $\tau'_{b,m_{s_j}}$ puts the possible pre-images of $y$ by $\phi\psi$  in a small interval, thus  we conclude as in Lemma \ref{tourc}. \qed\\

Note that Lemma \ref{tourab} gives only sufficient conditions, the same reasoning can produce many others.  It will not be used further, as Lemma \ref{tourc} is enough to prove

\begin{proposition}\label{mtours}
Let \begin{itemize} 
\item $\xi_n=\frac{1}{k_{n+2}}$ if the $n+1$-th multiplicative rule is  $I_m$ and $k_{n+1}\geq 2$,
\item $\xi_n=\frac{1}{3^lk_{n+2}...k_{n+l+1}}$ 
 if  the $n+1$-th  multiplicative rule is  $I_m$ with $k_{n+1}=1$ or $II_m$, and the next multiplicative rule   $I_m$ is the $n+l$-th, $l\geq 2$. \end{itemize}
Suppose $\sum\xi_n=+\infty$.
 Let $Z$ be the set of $y$ in  $Y_3$, such that $y$ is not in  $\tau'_{c,k}$ for all $k$ large enough. Then $\mu (Z)=0$ for the unique invariant measure $\mu$. \end{proposition}
{\bf Proof}\\
We fix a  multiplicative time $m_{n_0}$,  and for $n\geq n_0$ we define $Z_n$ to be the set of  $y$ which are not in  $\tau'_{c,k}$ for all $m_{n_0}\leq k\leq  m_n$, $n\geq n_0$, and $V_n$ such that $Z_n\setminus V_n=Z_{n+1}$. We have $Z_{n_0}=\tau'_{a,m_{n_0}} \cup \tau'_{b,m_{n_0}}$.\\

At each additive time $m_n\leq k < m_{n+1}$, the  new tower $\tau'_{c,k+1}$ is made with $\tau'_{c,k}$ stacked above one column of $\tau'_{a,k}$; $\tau'_{c,m_{n+1}-1}$ is made with $\tau'_{c,m_n}$ stacked above $k_{n+1}-1$ columns of $\tau'_{a,m_n}$; then, if the $n+1$-th multiplicative rule is $II_m$, $\tau'_{c,m_{n+1}}$ is made with $\tau'_{c,m_n}$ stacked above $k_{n+1}$ columns of $\tau'_{a,m_n}$;  if the $n+1$-th multiplicative rule is $I_m$, $\tau'_{c,m_{n+1}}$ is made with the last remaining column of $\tau'_{a,m_n}$. Then $V_n$ is made either with $k_{n+1}-1$ columns of $\tau'_{a,m_n}$ stacked above $\tau'_{c,m_n}$ plus
the last  column of $\tau'_{a,m_n}$, or with $k_{n+1}-1$ columns of $\tau'_{a,m_n}$ stacked above $\tau'_{c,m_n}$. In both cases, $V_n$ is a union of slices of  $\tau'_{a, m_{n+1}}$ and $\tau'_{c, m_{n+1}}$.\\

Assume, as is true for $n=n_0$, that  $Z_{n}$ is  a union of slices of $\tau'_{a, m_n}$ and $\tau'_{b, m_n}$; then $Z_n$ is also  a union of slices of  $\tau'_{a, m_{n+1}}$ and $\tau'_{b, m_{n+1}}$, and thus so is $Z_{n+1}$.\\
In all cases, $Z_{n+1}\cap \tau'_{a, m_{n+1}}$ is made with all $Z_{n}\cap \tau'_{b, m_{n}}$ and the intersection of $Z_{n}\cap \tau'_{a, m_{n}}$, with $k_{n+1}$  columns of $\tau'_{a, m_{n}}$ whose levels have measure $a_{m_{n+1}}$. $Z_{n+1}\cap \tau'_{b, m_{n+1}}$ is the intersection  of $Z_{n}\cap \tau'_{a, m_{n}}$ with one  column of $\tau'_{a, m_{n}}$ whose levels have measure $b_{m_{n+1}}$. Thus we have always, for $n\geq n_0+1$, $\mu(Z_{n}\cap \tau'_{a, m_{n}})\geq \mu(Z_{n}\cap \tau'_{b, m_{n}})$.\\

Figures 9 and 10 give a schematic view (note that the levels of the towers are not intervals, even when carried to $(X_9,T)$, see Figure 8 above) of what is used in the proof. The crossed parts form $V_n$, which has been deleted from $Z_n$ to get $Z_{n+1}$; the $\tau'_{c,m_n}$, crossed by dashed lines, have been deleted at an earlier stage.

\begin{figure}
\begin{center}
\begin{tikzpicture}[scale = 5]

\draw(0,.1)--(.6,.1);
\draw(1,.1)--(1.4,.1);
\draw(1.8,.1)--(2.1,.1);

\draw(0,.3)--(.6,.3);
\draw(1,.3)--(1.4,.3);
\draw(1.8,.3)--(2.1,.3);

\draw(0,.5)--(.6,.5);
\draw(1,.5)--(1.4,.5);

\draw(0,.7)--(.6,.7);
\draw(1,.7)--(1.4,.7);

\draw(0,.9)--(.6,.9);
\draw(1,.9)--(1.4,.9);

\draw(0,1.05)--(.6,1.05);

\draw(1,1.02)--(1.4,1.02);

\draw(0,.1)--(0,.5);
\draw(.6,.1)--(.6,.5);
\draw(1,.1)--(1,.5);
\draw(1.4,.1)--(1.4,.5);
\draw(1.8,.1)--(1.8,.3);
\draw(2.1,.1)--(2.1,.3);
\draw[dashed](0,.5)--(0,.7);
\draw[dashed](.6,.5)--(.6,.7);
\draw[dashed](1,.5)--(1,.7);
\draw[dashed](1.4,.5)--(1.4,.7);
\draw(0,.7)--(0,1.05);
\draw(.6,.7)--(.6,1.05);
\draw(1,.7)--(1,1.02);
\draw(1.4,.7)--(1.4,1.02);

\draw (.3,.25) node[below]{$\tau'_{a,m_n}$};

\draw (1.2,.25) node[below]{$\tau'_{a,m_n}$};

\draw (1.95,.25) node[below]{$\tau'_{a,m_n}$};

\draw (.3,.45) node[below]{$\tau'_{a,m_n}$};

\draw (1.2,.45) node[below]{$\tau'_{a,m_n}$};

\draw (.3,.85) node[below]{$\tau'_{a,m_n}$};

\draw (1.2,.85) node[below]{$\tau'_{a,m_n}$};

\draw (.3,1.03) node[below]{$\tau'_{b,m_n}$};

\draw (1.2,1) node[below]{$\tau'_{c,m_n}$};

\draw (.3,.1) node[below]{$\tau'_{a,m_{n+1}}$};

\draw (1.2,.1) node[below]{$\tau'_{b,m_{n+1}}$};

\draw (1.95,.1) node[below]{$\tau'_{c,m_{n+1}}$};

\draw(1.8,.1)--(2.1,.3);
\draw(1.8,.3)--(2.1,.1);
\draw(1,.3)--(1.4,.9);
\draw(1,.9)--(1.4,.3);
\draw[dashed](1,.9)--(1.4,1.02);
\draw[dashed](1,1.02)--(1.4,.9);

\end{tikzpicture}\\
\caption{Cutting and stacking $I_m$}

\end{center}

\end{figure}
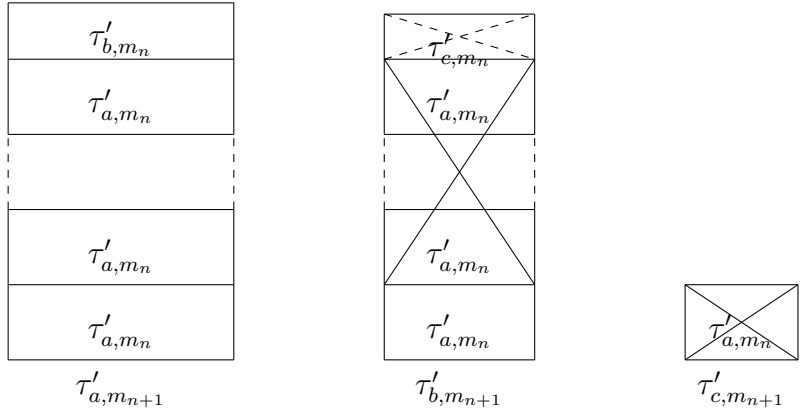

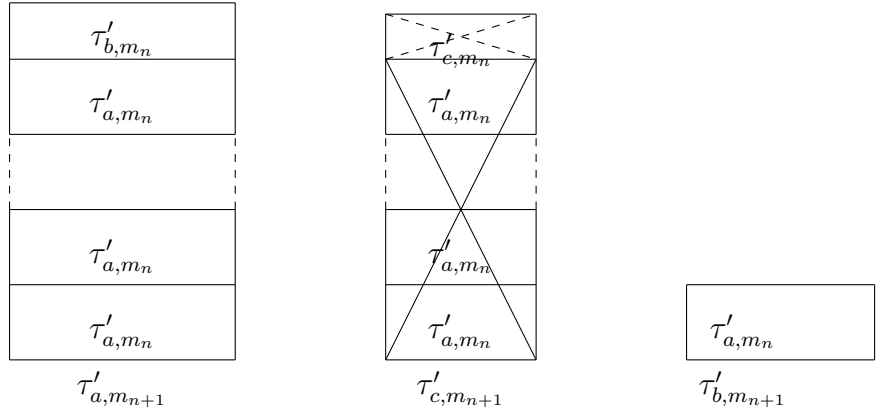
\begin{figure}[h] \label{}
\begin{center}
\begin{tikzpicture}[scale = 5]

\draw(0,.1)--(.6,.1);
\draw(1,.1)--(1.4,.1);
\draw(1.8,.1)--(2.3,.1);

\draw(0,.3)--(.6,.3);
\draw(1,.3)--(1.4,.3);
\draw(1.8,.3)--(2.3,.3);

\draw(0,.5)--(.6,.5);
\draw(1,.5)--(1.4,.5);

\draw(0,.7)--(.6,.7);
\draw(1,.7)--(1.4,.7);

\draw(0,.9)--(.6,.9);
\draw(1,.9)--(1.4,.9);

\draw(0,1.05)--(.6,1.05);

\draw(1,1.02)--(1.4,1.02);

\draw(0,.1)--(0,.5);
\draw(.6,.1)--(.6,.5);
\draw(1,.1)--(1,.5);
\draw(1.4,.1)--(1.4,.5);
\draw(1.8,.1)--(1.8,.3);
\draw(2.3,.1)--(2.3,.3);
\draw[dashed](0,.5)--(0,.7);
\draw[dashed](.6,.5)--(.6,.7);
\draw[dashed](1,.5)--(1,.7);
\draw[dashed](1.4,.5)--(1.4,.7);
\draw(0,.7)--(0,1.05);
\draw(.6,.7)--(.6,1.05);
\draw(1,.7)--(1,1.02);
\draw(1.4,.7)--(1.4,1.02);

\draw (.3,.25) node[below]{$\tau'_{a,m_n}$};

\draw (1.2,.25) node[below]{$\tau'_{a,m_n}$};

\draw (1.95,.25) node[below]{$\tau'_{a,m_n}$};

\draw (.3,.45) node[below]{$\tau'_{a,m_n}$};

\draw (1.2,.45) node[below]{$\tau'_{a,m_n}$};

\draw (.3,.85) node[below]{$\tau'_{a,m_n}$};

\draw (1.2,.85) node[below]{$\tau'_{a,m_n}$};

\draw (.3,1.03) node[below]{$\tau'_{b,m_n}$};

\draw (1.2,1) node[below]{$\tau'_{c,m_n}$};

\draw (.3,.1) node[below]{$\tau'_{a,m_{n+1}}$};

\draw (1.2,.1) node[below]{$\tau'_{c,m_{n+1}}$};

\draw (1.95,.1) node[below]{$\tau'_{b,m_{n+1}}$};

\draw(1,.1)--(1.4,.9);
\draw(1,.9)--(1.4,.1);
\draw[dashed](1,.9)--(1.4,1.02);
\draw[dashed](1,1.02)--(1.4,.9);

\end{tikzpicture}\\
\caption{Cutting and stacking $II_m$}

\end{center}

\end{figure}

We want now to estimate the measure of $V_n$.\\

We suppose first that the $n+1$-th multiplicative rule is  $I_m$. If $k_{n+1}\geq 2$, $V_n$ is a slice of $\tau'_{b,m_{n+1}}$ of height $(k_{n+1}-1)h_{a,m_n}$. If $k_{n+1}=1$, $V_n$ is $\tau'_{c,m_{n+1}}$.

Suppose $k_{n+1}\geq 2$. Then we need to  estimate $\mu (\tau'_{b,m_{n+1}})$; we notice that $\mu (\tau'_{a,m_{n+2}})\geq \frac{1}{5}$, because this tower is wider than the two others, and at least half as high by the estimate at the end of Section \ref{pqm}. $\tau'_{b,m_{n+1}}$, is a slice of $\tau'_{a,m_{n+2}}$ of height 
$h_{b,m_{n+1}}$, while $h_{a,m_{n+2}}=k_{n+2}h_{a,m_{n+1}}+h_{b,m_{n+1}}$. From $h_{b,m_{n+1}}=k_{n+1}h_{a,m_n}+h_{c,m_n}$, $h_{a,m_{n+1}}=k_{n+1}h_{a,m_n}+h_{b,m_n}$,
we get $h_{b,m_{n+1}}\geq \frac{k_{n+1}}{k_{n+1}+2}h_{a,m_{n+1}}\geq \frac{1}{2}h_{a,m_{n+1}}$, and $\mu (\tau'_{b,m_{n+1}})\geq \frac{1}{10(k_{n+2}+2)}$. Now $V_n$ is a slice of $\tau'_{b,m_{n+1}}$ of relative height at least $\frac{k_{n+1}-1}{k_{n+1}+2}\geq \frac{1}{3}$, 
and we get $\mu (V_n)\geq \frac{1}{30(k_{n+2}+2)}$.

If $k_{n+1}=1$, we take first $l=2$: the $n+2$-th multiplicative rule is also $I_m$. Then $\tau'_{c,m_{n+1}}$, is a slice of $\tau'_{b,m_{n+2}}$ of height 
$h_{c,m_{n+1}}$, while $h_{b,m_{n+2}}=k_{n+2}h_{a,m_{n+1}}+h_{c,m_{n+1}}$. We have $h_{a,m_{n+1}}=h_{a,m_n}+h_{b,m_n}$, $h_{c,m_{n+1}}=h_{a,m_n}$, thus $h_{c,m_{n+1}}\geq \frac{1}{3}h_{a,m_{n+1}}$, thus we get 
$\mu (\tau'_{c,m_{n+1}})\geq \frac{1}{3(k_{n+2}+2)}\mu (\tau'_{b,m_{n+2}})$. Then $\mu (\tau'_{b,m_{n+2}})$ is estimated just as $\mu (\tau'_{b,m_{n+1}})$ in the case above, with the only difference  that $k_{n+2}$ may be equal to one: we get it is at least
$\frac{1}{15(k_{n+3}+2)}$, and thus $\mu(V_n)\geq \frac{1}{45(k_{n+2}+2)(k_{n+3}+2)}$.

For larger values of $l$ we iterate this method, looking at $\tau'_{c,m_{n+1}}$, inside ... inside $\tau'_{c,m_{n+l-1}}$, inside $\tau'_{b,m_{n+l}}$,inside $\tau'_{a,m_{n+l+1}}$, Estimating the measures gives us first factors
$\frac{1}{k_{n+2}+2}$ ... $\frac{1}{k_{n+l+1}+2}$, but depend also on the comparison of successive heights of towers, which brings  factors $\frac13$.

If the $n+1$-th multiplicative rule is $II_m$, $V_n$ is a slice of $\tau'_{c,m_{n+1}}$ of height at least $\frac{1}{3}h_{c,m_{n+1}}$ and we estimate its measure in the same way.\\

In all cases $V_n$ is a columns- of $\tau'_{a,m_n}$ while $Z_n\cap \tau'_{a,m_n}$ is a slice of $\tau'_{a,m_n}$, and for any column $\Lambda$ and slice $\Lambda'$ of the same tower we have $\mu(\Lambda\cup \Lambda')=\mu(\Lambda)\mu(\Lambda')$. Thus $\mu(V_n\cap Z_n)\geq \mu(Z_n\cap\tau'_{a,m_n}\cap V_n )=\mu(V_n)\mu(Z_n\cap\tau'_{a,m_n})\geq \frac{1}{2}\mu(V_n)\mu(Z_n)$, 
and $\mu(Z_{n+1})\leq \mu(Z_n)(1-\frac{1}{2}\mu(V_n))\leq \mu(Z_n)(1-K\xi_n)$ for some constant $K$, and we conclude by a Borel-Cantelli argument, namely 
$\mu(Z_{n+1})\leq \mu(Z_{n_0})\prod_{n>n_0}(1-K\xi_n)$, thus $\mu(Z)=0$ because of our hypothesis. \qed\\

We turn now to the isomorphism problem: as $E_3$ is nonempty, the best we can hope is to replace the semi-conjugacies in Section \ref{nno} by measure-theoretic isomorphisms.

\begin{theorem}\label{imp} Under the hypothesis of Proposition \ref{mtours}, an AR9 or AR6 symbolic system or interval exchange is uniquely ergodic and measure-theoretically isomorphic to its AR3 coding. \end{theorem}
{\bf Proof}\\
Then, by Proposition \ref{mtours} and Lemma \ref{tourc} $\phi$ is invertible almost everywhere. Thus $\phi$ provides a measure-theoretic isomorphism between $(Y_3,S,\mu)$ and $(Y_9,S,\mu')$ for any normalized invariant measure $\mu'$. Such an invariant measure $\mu'$ can be defined also on $(X_9,T)$  as $\psi$ is invertible almost everywhere, and  $\psi$ provides a measure-theoretic isomorphism between $(X_9,T,\mu')$ and $(Y_9,S,\mu')$. In particular, any such measure $\mu'$ has to be ergodic, hence the unique ergodicity. The results extend then to the  intermediate coding $(Y_6,S,\mu')$ and to its geometric model  $(X_6,T,\mu')$.  \qed\\

\begin{definition}\label{aa}  As in \cite{bst}, we consider measures  on  all infinite sequences of symbols  $I$, $II$, $III$  and take any shift invariant ergodic probability measure  $\nu$ which assigns positive measure to each cylinder;  by identifying an AR3, AR6,  or AR9 system with its defining sequence $(r_n)$, we can define $\nu$ on the set of all AR3, AR6, or  AR9 systems.\end{definition}

In particular, one of these  measures  coincides with  the measure of maximal entropy for the suspension flow of the Rauzy gasket built in \cite{ahs0}, see also \cite{ahs1}. 
 
\begin{proposition}\label{pt} The hypothesis of Proposition \ref{mtours} is satisfied by $\nu$-almost every AR3, AR6, or AR9 system.\end{proposition}
{\bf Proof}\\  This hypothesis is satisfied in particular if for infinitely many $n$ we have $k_{n+1}=2$ and $k_{n+2}=1$, which is satisfied in particular if for infinitely many $p$ we have $r_p=I$, $r_{p+1}=r_{p+2}=III$, $r_{p+3}=r_{p+4}=I$. As this cylinder has positive measure and $\nu$ is ergodic, this is true for $\nu$-almost every sequence $(r_n)$.\qed\\

This completes the proof of Theorem \ref{main} above. But the sufficient condition in Proposition \ref{mtours} gives also the isomorphism (and unique ergodicity) for many explicit examples; while the  first set of values of $\xi_n$ is enough to prove Proposition \ref{pt} above, with the help of the second set of values  we can prove the following. 

\begin{proposition}\label{bqp} The hypothesis of Proposition \ref{mtours} is satisfied by all Arnoux-Rauzy systems where the $k_n$ are bounded (in \cite{bst} these are said to have {\em bounded weak partial quotients}). \end{proposition}
{\bf Proof}\\ If the $n+2$-th multiplicative rule is  $I_m$, then $\xi_n$ is either $\frac{1}{k_{n+2}}$ or $\frac{1}{3^2k_{n+2}k_{n+3}}$. As there are infinitely many rules $I_m$, we get infinitely many $n$ for which $\xi_n\geq \frac{1}{9K_0^2}$ if all the $k_n$ are bounded by $K_0$. \qed \\

This completes the proof of Corollary \ref{cp1} above; then  Corollary \ref{cp2} is proved by using the measure-theoretic  isomorphism between the Tribonacci AR3 and a rotation of the $2$-torus \cite{rau} and the fact that such a rotation is always rigid.

\subsection{Non unique ergodicity}\label{snue}

\begin{theorem}\label{nue} If $\sum_{n=1}^{+\infty} \frac{1}{k_n}<+\infty$, each corresponding AR9 or AR6 symbolic system or interval exchange is not uniquely ergodic; it has two ergodic invariant measures; it is measure-theoretically isomorphic to its AR3 coding if and only if  it is equipped with an ergodic measure. \end{theorem}
{\bf Proof}\\
Let $\mu'$ be any normalized invariant measure on $(Y_9,S)$. We first show that at multiplicative times all towers have very small measure except $\tau'_{1,m_n}$ and $\tau'_{4,m_n}$.

Indeed, from the multiplicative rules of Section \ref{pqm} we get that $\tau'_{b,m_n}$ is a slice of $\tau'_{a,m_{n+1}}$ of height $h_{b,m_n}$, hence $\mu (\tau'_{b,m_n})\leq \frac{2}{k_{n+1}-1}$, while $\tau'_{c,m_n}$ is a slice of either  $\tau'_{b,m_{n+1}}$ or  $\tau'_{c,m_{n+1}}$, of height $h_{c,m_n}$, hence $\mu (\tau'_{c,m_n})\leq \frac{2}{k_{n+1}-1}$; and $\mu'(\tau'_{i,m_n})\leq \mu'(\psi\tau_{b,m_n})=\mu(\tau'_{b,m_n})$ for $i=5,6,7$, $\mu'(\tau'_{i,m_n})\leq \mu'(\psi\tau_{c,m_n})=\mu(\tau'_{c,m_n})$ for $i=8,9$. \\
Now, from the multiplicative rules at the end of Section \ref{ccr} we get that $\tau'_{2,m_n}$ is either $\tau'_{6,m_{n+1}}$ or the union  of $\tau'_{8,m_{n+1}}$ with a slice of 
$\tau'_{7,m_{n+1}}$, thus $\mu'(\tau'_{2,m_n})\leq \frac{4}{k_{n+2}-1}$. Finally $\tau'_{3,m_n}$ is either the union  of $\tau'_{9,m_{n+1}}$ with a slice of 
$\tau'_{1,m_{n+1}}$ of relative height at most $\frac{1}{k_{n+1}-1}$, or the union  of $\tau'_{5,m_{n+1}}$ with a slice of 
$\tau'_{4,m_{n+1}}$ of relative height at most $\frac{1}{k_{n+1}-1}$: in both cases $\mu'(\tau'_{3,m_n})\leq \frac{3}{k_{n+1}-1}$.\\

Thus, the condition $\sum_{n=1}^{+\infty} \frac{1}{k_n}<+\infty$ implies that for any invariant measure $\mu'$, the system $(Y_9,S,\mu')$ is such that $\mu'$-almost every point $y$ in $Y_9$ is determined by the sequences $\iota"(y,k)$,  $\eta(y,k)$ such that $y$ is in level $\eta(y,k)$ of the tower $\tau'_{\iota"(y,k), k}$, $\iota"(y,k)\in \{1,4\}$. We say that $(Y_9,S,\mu')$ is {\em generated} by two sequences of towers, and such a measure-theoretic system is said to be a system of {\em rank} (at most) two; by a classical result for which we refer the reader to \cite{fr1}, $(Y_9,S)$, which is of rank at most two for any invariant measure, has as at most two ergodic invariant measures.\\

At multiplicative times, we define  recursively  $(\tau'_{{\bar 1},m_n},\tau'_{{\bar 4},m_n})=e^l(\tau'_{{1},m_n},\tau'_{{4},m_n})$ if $l$ is the total number of rules $I_m$ (strictly) before the $n$-th multiplicative  rule and $e$ is the exchange. Then for each $n$, $\tau'_{{\bar 1},m_n}$ makes all but a very small part of $\tau'_{{\bar 1},m_{n+1}}$, $\tau'_{{\bar 4},m_n}$ makes all but a very small part of $\tau'_{{\bar 4},m_{n+1}}$, and all the other  $\tau'_{i,m_n}$, $i\neq 1,4$ have very small measure.\\

We define a new  symbolic  system $({\bar X}, {\bar T}, {\bar \mu})$ on the alphabet $\{a,s\}$ by the words $D_0=a$,
 $D_{n+1}=s^{h_{a,m_n}}D_n^{k_{n+1}-1}s^{h_{b,m_n}}$. By a standard argument, see \cite{fr1}, we can build 
 towers $\bar\tau'_n$ in $\bar X$, $\bar\tau'_{n+1}$ being obtained from $\bar\tau'_n$ by cutting  it into $k_{n+1}-1$ equal columns, stacking them above each other, stacking below them $h_{a,m_n}$ new levels called spacers, and stacking above them $h_{b,m_n}$ new levels called spacers; almost every point $x$ in $\bar X$ is determined by the sequence  $\eta'(x,n)$ such that $y$ is in level $\eta'(x,n)$ of the tower $\bar\tau'_n$.  $({\bar X}, {\bar T}, {\bar \mu})$ is a {\em system of rank one}, as it can be generated by a single family of towers.
 
As is explained in more details in \cite{afp}, we can build an application $\phi_1$ from $\bar X$ to $Y_9$ 
 by sending the $j$-th level of 
the tower $\bar\tau'_{n}$ to the $j$-th level of the  tower $\tau'_{{\bar 1},m_n}$: it is
consistent,  defined almost
everywhere and one-to-one. By taking the image of $\bar \mu$ by $\phi_1$,  we build
 a measure-theoretic 
isomorphism
between the rank one system  $({\bar X}, {\bar T}, {\bar \mu})$  and 
 $(Y_9, S)$ equipped with some invariant probability measure $\mu_{1}$; $\mu_{1}$ is ergodic as $\bar \mu $ is. We do the same for 
 another application $\phi_4$, which sends the $j$-th level of 
$\bar\tau'_n$ to the $j$-th level of $\tau'_{{\bar 4},m_n}$.
 defining an ergodic $\mu_{4}$. Now, $\mu_1(\tau'_{{\bar 1},m_n})$ and $\mu_4(\tau'_{{\bar 4},m_n})$ are close to $1$, $\mu_1(\tau'_{{\bar 4},m_n})$ and $\mu_4(\tau'_{{\bar 1},m_n})$ are close to $0$ for $n$ large enough, thus there exists $n$ for which $\mu_1(\tau'_{{\bar 1},m_n})\neq \mu_4(\tau'_{{
 \bar 1},m_n})$, thus $\mu_1\neq \mu_4$ on $(Y_9,S)$. \\
 
 The results extend immediately to $(X_9,T)$, and to the AR6 systems, to which we carry $\mu_1$ and $\mu_4$. \\

 Now, the AR3 coding $(Y_3,S,\mu)$ is also a  system of rank one, generated by the towers $\tau'_{a,m_n}$. These towers are built in the same way as the $\bar\tau'_n$, as replacing a small part of $\tau'_{a,m_n}$ by spacers does not change the system, thus as in \cite{afp} $(Y_3,S,,\mu)$ is measure-theoretically isomorphic to $({\bar X}, {\bar T}, {\bar \mu})$, thus to both $(Y_9,S,\mu_1)$ and  $(Y_9,S,\mu_4)$; but it cannot be  measure-theoretically isomorphic to a non-ergodic $(Y_9,S,\mu')$. And the same reasoning holds for the others AR9 or AR6 systems considered.
 \qed\\

Note that  in  the only family of counter-examples we have, the two-point extension of Proposition \ref{E123} is rather degenerate, being ergodic only when the measure is concentrated on one copy of the factor.

\section{Weak mixing}\label{wm}

 \begin{definition}\label{dw}
If
$(X', U, \mu_0)$ is a finite
measure-preserving dynamical system, a real number $0\leq \theta <1$ is a {\em measurable eigenvalue} 
(denoted additively) if there
exists a non-constant $f$ in ${\mathcal L}^1(X', \GR/\GZ )$  such
that
$f\circ U= f+ \theta$ (in ${\mathcal L}^1(X', \GR/\GZ )$); $f$ is
then
 an {\em eigenfunction} for the eigenvalue $\theta$.
 
 As
constants are not eigenfunctions,
$\theta =0$ is not an eigenvalue if $U$ is  ergodic.

 $(X', U, \mu_0)$ is {\em
weakly mixing}
if it has no measurable eigenvalue.
\end{definition}

The existence of weak mixing for AR3 systems, proved in \cite{cfme}, came as a surprise;  this existence persists for AR9 (and AR6) systems, because under the hypothesis $\sum_{n=1}^{+\infty} \frac{1}{k_n}<+\infty$, by Theorem \ref{nue} above the AR9  or AR6 system equipped with one of its ergodic measures is isomorphic to its AR3 coding, while by Theorem 2 of \cite{cfme} this  AR3 system is weakly mixing. The sufficient condition given in \cite{cfme} for weak mixing of AR3 systems is weaker than the condition $\sum_{n=1}^{+\infty} \frac{1}{k_n}<+\infty$:  we shall show now that under this sufficient  conditions the AR9 systems are also weakly mixing, for any ergodic invariant measure. But indeed this raises more questions than gives answers, as we shall see in the discussion below.

\begin{proposition}\label{twm} An ergodic AR9  or AR6 system  is weakly mixing if \begin{itemize}
\item $k_{n_i+2}$ is unbounded,
\item $$\sum_{i=1}^{+\infty}\frac{1}{k_{n_i+1}}<+\infty,$$
\item  $$\sum_{i=1}^{+\infty}\frac{1}{k_{n_i}}<+\infty,$$ \end{itemize}
where the $n_i$ are the $n\geq 1$ for which the $n$-th multiplicative rule is $I_m$. \end{proposition}

{\bf Proof}\\
The only difference between the present proof and the proof in \cite{cfme} is in the beginning. Namely, to prove Proposition 10 of \cite{cfme}, we use the fact that when we move by $S^{h_{a,m_n}}$ inside a substantial slice of $\tau'_{a,m_{n+1}}$, we arrive at the same level in $\tau'_{a,m_n}$; here we need the stronger result that for all $i=1,2,3,4$, when we move by $S^{h_{a,m_n}}$ inside a substantial slice of $\tau'_{i,m_{n+1}}$, we arrive at the same level in some $\tau'_{j,m_n}$. This in turn involves some technical difficulties when $k_{n+1}$ is small, obligeing us to use our hypotheses on the $k_n$ at that stage, which was not recessary in \cite{cfme}. Thus Proposition 10 of \cite{cfme} is replaced by 

  \begin{lemma}\label{chnad} If  $\theta$ is a measurable eigenvalue for an  AR9  symbolic system $(Y_9,S,\mu')$ satisfying the hypotheses of Proposition \ref{twm}, 
  $k_{n+1}\ab\ab h_{a,m_n} \theta \ab\ab \to 0$ when $n\to+\infty$,
  where $\ab\ab \quad\ab\ab $ denotes the distance to the nearest
  integer.\end{lemma}
  {\bf Proof}\\
  Let $f$ be an eigenfunction for the eigenvalue $\theta$; for each $\varepsilon>0$ there exists $N(\varepsilon )$ such that for all $n>N(\varepsilon)$ there exists $f_n$, which  satisfies $\int||f-f_n||d\mu <\varepsilon$
and is constant on each level of each tower $\tau'_{i,m_{n-2}}$, $\tau'_{i,m_{n-1}}$, and $\tau'_{i,m_n}$, $i=1,...9$.\\

 Suppose  first $k_{n+1}\geq 3$. Let $j$ be any integer with $0\leq j\leq\left[\frac{k_{n+1}-1}2\right]$.
  
  Suppose for example the $n+1$-th multiplicative rule is  $I_m$; we have the concatenation rule $1_{m_{n+1}}=3_{m_n}4_{m_n}^{k_{n+1}-1}5_{m_n}$.
Let  $\tau''_n$ be the slice of $\tau'_{1,m_{n+1}}$ consisting of levels  from $h_{a,m_n}$ to $h_{a,m_n}+{[\frac{k_{n+1}-1}2]h_{a,m_n}-1}$;  it has relative height at least $\frac15$.

By construction, for any point $x$ in $\tau''_n$, $S^{jh_{a,m_n}}x$ is in the tower $\tau'_{1,m_{n+1}}$, and in the same level of the tower $\tau'_{4,m_n}$ as $x$.
Thus for $\mu'$-almost every $x\in \tau''_{n}$, $f_{n}(S^{jh_{a,m_n}}x)=f_{n}(x)$ while
$f(S^{jh_{a,m_n}}x)=\theta jh_{a,m_n}+f(x)$;
we have
\[
\int_{\tau''_{n}}\ab\ab f_{n}\circ S^{jh_{a,m_n}}-j\theta h_{a,m_n}-f_{n}\ab\ab d\mu'=
\int_{\tau''_{n}}\ab\ab j\theta h_{a,m_n}\ab\ab d\mu' = \ab
\ab j\theta h_{a,m_n}\ab\ab \mu (\tau''_n)
\]
and
\[\int_{\tau''_{n}}\ab\ab f_{n}\circ S^{jh_{a,m_n}}-j\theta h_{a,m_n}-f_{n}\ab\ab d\mu' \leq
\int_{\tau''_{n}}\ab\ab f_{n}\circ S^{jh_{a,m_n}}-f\circ S^{jh_{a,m_n}}\ab\ab d\mu'+
\int_{\tau''_{n}}\ab\ab f_{n}-f\ab\ab d\mu'
< 2\varepsilon.
\]

Thus we get
$\ab\ab j\theta h_{a,m_n}\ab\ab \mu'(\tau'_{1,m_{n+1}})<10 \varepsilon$, for $n>N(\varepsilon)$ and any integer $0\leq j\leq \left[\frac{k_{n+1}-1}2\right]$.

The same result holds when the $n+1$-th multiplicative rule is $II_m$, with concatenation rule $1_{m_{n+1}}=1_{m_n}^{k_{n+1}}7_{m_n}$: just $\tau'_{1,m_n}$ replaces 
$\tau'_{4,m_n}$. And the same construction, mutatis mutandis, works with $\tau'_{1,m_{n+1}}$ replaced by $\tau'_{i,m_{n+1}}$, $i=2,3,4$. Summing the four inequalities and taking into account that $\sum_{i=1}^4\mu'(\tau'_{i,m_{n+1}})=\mu (\tau'_{a,m_{n+1}})\geq \frac15$, we get
$\ab\ab j\theta h_{a,m_n}\ab\ab <50 \varepsilon$ for $0\leq j\leq \left[\frac{k_{n+1}-1}2\right]$, hence
$\ab\ab j\theta h_{a,m_n}\ab\ab <200 \varepsilon$ for $0\leq j\leq k_{n+1}$.

We continue exactly as in \cite{cfme}. Let $\varepsilon <\frac{1}{1000}$, and suppose $\ab\ab k_{n+1}\theta h_{a,m_n}\ab\ab \neq k_{n+1}\ab\ab \theta h_{a,m_n}\ab\ab$: let $i$ be the smallest $0\leq j\leq k_{n+1}$ such that $\ab\ab j\theta h_{a,m_n}\ab\ab \neq j\ab\ab \theta h_{a,m_n}\ab\ab$, then $i\geq 2$ and
$\ab\ab (i-1)\theta h_{a,m_n}\ab\ab =(i-1)\ab\ab \theta h_{a,m_n}\ab\ab$, thus
$i\ab\ab \theta h_{a,m_n}\ab\ab =(i-1)\ab\ab \theta h_{a,m_n}\ab\ab+\ab\ab \theta h_{a,m_n}\ab\ab=\ab\ab (i-1)\theta h_{a,m_n}\ab\ab +\ab\ab \theta h_{a,m_n}\ab\ab<400\varepsilon<\frac{1}{2}$ thus
$\ab\ab i\theta h_{a,m_n}\ab\ab =\ab\ab (i\ab\ab \theta h_{a,m_n}\ab\ab)\ab\ab =i \ab\ab \theta h_{a,m_n}\ab\ab$, contradiction. Thus we get $k_{n+1}\ab\ab \theta h_{n-1}\ab\ab< 200 \varepsilon$  for $n>N(\varepsilon)$.\\

Suppose now $k_{n+1}=2$; then, except maybe for a finite number of values of $n$, the hypotheses imply that the $n$-the multiplicative rule is  $II_m$. Note also that we need only to prove  $\ab\ab \theta h_{a,m_n}\ab\ab <C \varepsilon$. For concatenation rules such as
$1_{m_{n+1}}=1_{m_n}^27_{m_n}$, we see that $S$ iterated by the length of $1_{m_n}$, namely $h_{a,m_n}$, sends to itself each level of 
$\tau'_{1,m_n}$ if we start from the first slice $\tau'_{1,m_n}$ in $\tau'_{1,m_{n+1}}$, whose height is comparable (by some constant) to the height of $\tau'_{1,m_{n+1}}$, thus we can write the reasoning which leads to $\ab\ab \theta h_{a,m_n}\ab\ab \mu'(\tau'_{1,m_{n+1}})<C\varepsilon$. \\
If there is no square in the concatenation rule, its right member is $3_{m_n}4_{m_n}5_{m_n}$, which is equal to
$4_{m_{n-1}}^{k_n}5_{m_{n-1}}3_{m_{n-1}}4_{m_{n-1}}^{k_n-1}5_{m_{n-1}}3_{m_{n-1}}$; if $k_n\geq 2$ we iterate $S$ by the length of $4_{m_{n-1}}^{k_n}5_{m_{n-1}}$, which is $h_{a,m_n}$,
starting from the  $4_{m_{n-1}}^{k_n-1}$, at the end of $4_{m_{n-1}}^{k_n}$; if $k_n =1$ we use the length  of $3_{m_{n-1}}5_{m_{n-1}}$, which is $h_{a,m_n}$,
starting from the  first $3_{m_{n-1}}$. In both cases, the iteration of $S$ by the chosen quantity will send levels of some $\tau'_{i,m_{n-1}}$ to themselves, thus our choice of $f_n$ allows to write the usual reasoning, and to complete the case $k_{n+1}=2$.\\

Suppose $k_{n+1}=k_n=1$.Then, again for $n$ large enough, the $n$-th and $n-1$-th multiplicative rules are  $II_m$. The concatenations we look for are $1_{m_n}7_{m_n}=1_{m_{n-1}}7_{m_{n-1}}1_{m_{n-1}}$, $4_{m_n}5_{m_n}=3_{m_{n-1}}5_{m_{n-1}}3_{m_{n-1}}$, $3_{m_n}5_{m_n}=4_{m_{n-1}}5_{m_{n-1}}3_{m_{n-1}}$,$4_{m_n}6_{m_n}=3_{m_{n-1}}5_{m_{n-1}}2_{m_{n-1}}$, In the first one, $h_{a,m_n}$, which is the length of $1_{m_{n-1}}7_{m_{n-1}}$, can be used to iterate $S$ starting from the first $1_{m_{n-1}}$, and similarly in the second one. The last ones are  equal to 
$3_{m_{n-2}}4_{m_{n-2}}^{k_{n-1}-1}5_{m_{n-2}}3_{m_{n-2}}4_{m_{n-2}}^{k_{n-1}}5_{m_{n-2}}$ and 
$4_{m_{n-2}}^{k_{n-1}-}5_{m_{n-2}}3_{m_{n-2}}4_{m_{n-2}}^{k_{n-1}}6_{m_{n-2}}$: in both cases  
$h_{a,m_n}$, which is the length of $4_{m_{n-2}}^{k_{n-1}-}5_{m_{n-2}}3_{m_{n-2}}$, can be used to iterate $S$ starting from a sizeable slice of the tower. \\

Suppose $k_{n+1}=1$ but $k_n\geq 2$. Then $h_{a,m_n}=h_{a,m_{n+1}}-h_{a,m_{n-1}}$. We have $\ab\ab \theta h_{a,m_{n-1}}\ab\ab <C \varepsilon$ because $k_n\geq 2$, and $\ab\ab \theta h_{a,m_{n+1}}\ab\ab <C \varepsilon$ either because $k_{n+2}\geq 2$ or because 
$k_{n+1}=k_{n+2}=1$, thus we conclude.\qed\\

Then the (nontrivial) Sections 3 and 4 of \cite{cfme} prove that, under the hypotheses of Proposition \ref{twm}, the condition $k_{n+1}\ab\ab h_{a,m_n} \theta \ab\ab \to 0$ gives no possible $\theta$ except $\theta=0$, which is excluded because of the ergodicity of the system. The same reasoning applies to the other AR9 or AR6 systems. \qed
\\

We do not know whether this sufficient condition gives interesting new examples; it might help  to find a weakly mixing AR9 system for which $\mu(E_1)=1$ in the AR3 coding, but this we were not able to achieve. Indeed, starting from Lemma \ref{tourc} as in Section \ref{simp}, we are able to build such AR9 systems under the condition  $\sum_{i=1}^{+\infty}\frac{1}{k_{n_i+1}}=+\infty.$ while $\sum_{i=1}^{+\infty}\frac{1}{k_{n_i}}$ may be finite; we could also get these conditions by  starting from Lemma \ref{tourab} and  imitating the proof of Proposition \ref{mtours}; this  falls short of being compatible with the conditions of Proposition \ref{twm}. Indeed, we conjecture that these conditions are not compatible with $\mu(E_1)=1$, and not even with unique ergodicity; whether these conditions are necessary for weak mixing is not known either. It would be also very interesting to find a uniquely ergodic weakly mixing AR9, or a weakly mixing AR9 which is not isomorphic to its AR3 coding.

\end{document}